\input amstex
\documentstyle{amsppt}
\topmatter
\title 
Automorphisms of Finite Order on Rational Surfaces \\
\\({with an Appendix by I. Dolgachev})
\endtitle
\author 
D. -Q. Zhang
\endauthor 
\endtopmatter
\NoBlackBoxes
%
%
%
%
%
\par \vskip 1pc \hskip 2pc
Department of Mathematics, National University of Singapore, Singapore
\par \vskip 1pc \hskip 8.5pc
e-mail : matzdq$\@$math.nus.edu.sg

\par \vskip 3pc 
We classify minimal pairs $(X, G)$
for smooth rational projective surface $X$ and
finite group $G$ of automorphisms
on $X$. We also determine the fixed locus $X^G$
and the quotient surface $Y = X/G$ as well as
the fundamental group of the smooth part of $Y$.
The realization of each pair is
included.  Mori's extremal ray theory and recent results
of Alexeev and also Ambro on
the existence of good anti-canonical divisors are used.


\par \vskip 3pc
{\bf Introduction}

\par \vskip 1pc
More than one hundred years ago, Kantor had written a book on finite
birational automorphism groups of rational surfaces.
In the sixties to eighties, Manin, Iskovskih, Gizatullin also thoroughly 
studied  $G$-rational surfaces defined over non-closed fields.
One aim of them is to reduce to $G$-minimal surfaces.
In [Giz], $G$-pseudoprojective rational surfaces, which are not
$G$-projective surfaces, are shown to be
relative minimal elliptic surfaces; the same paper also shows that
not every  $G$-rational surface is  $G$-pseudo-projective.
B. Segre [Seg] did, among many other things, the classification
of $Aut X$  for cubic surfaces $X$ (see also [Ho2]).  

\par
Recently, $Aut X$  has also been classified again for quartic
del Pezzo surface in [Ho1].  In [Koit], automorphism groups of rational
surfaces obtained by blowing up very general points in  $\bold P^2$
are completely classified.  
It is very desirable to test the modern machineries on the
old subject and obtain a simpler proof at the same time.

\par
In this note, we work over the complex numbers field  $\bold C$
and consider pairs of  $(X, G)$  of an arbitrary 
smooth rational projective surface  $X$  with a fixed finite group  
$G$  acting on it.  To simplify the arguments, we assume also that  $G$
is cyclic of prime order.  We believe that the general case could be handled
similarly.  Indeed, our last theorem deals with arbitrary  $G$,
where we reduce to either  $G$-stable conic
fibration or del Pezzo case (see Remark 5).

\par
Actually, this note is inspired by Bayle-Beauville's recent simple new 
classification
of birational involutions of rational surfaces [BB].  As there, we also adopt the
latest Mori-theory [Mor].  Though the theory has been developed along the course
of classification of higher dimensional varieties (dimension at least 3),
we will see how useful it is even for surfaces.
First, it will help us to reduce to a $G$-minimal surface
very quickly, which has either a $G$-stable conic-fibration (Mori-fibration),
or a Picard number one quotient surface.  The first case is easy to treat.

\par
For the second case, we have two approaches.  The top down approach is
based on known informations on Weyl groups  $W(E_n)$  of lattice  $E_n$
where we apply Manin's results in [Man2]; see also [Re1].
For the bottom up approach (more geometric), 
we will study the quotient surfaces;
this approach is normally more difficult; 
to do so, we apply results of Alexeev and Ambro about the existence of
a good member in the anti-canonical system [Alex], [Am]; implicitly 
we are also using
Fujita's theory of polarized varieties : $\Delta$ genus zero case
[Fuj]; this way, we avoid referring to the classification list
of automorphism groups of del Pezzo surfaces  $X$; such a list is available
if  $K_X^2$  is bigger. 

\par
It turns out that all pairs  $(X, G)$  with  $G$  cyclic of prime order  $p$, 
except the last 3 rows in Table 1 ($p = 5$), have minimal models  
$(X_{\min}, G)$, via a  $G$-equivariant birational morphism 
(only smooth blow-downs of  $G$-stable divisors
but no blow-ups), such that at least one of  $X$  and  $Y = X/G$  is 
a minimal rational surface (i.e., $\bold P^2$  or Hirzebruch surfaces  $\bold F_e$,
$e \ne 1$) or the projective cone  $\overline{\bold F}_3$ ($p = 3$).

\par
To be precise, denote by  $\mu_p$  the multiplicative group of prime order  $p$.
By writing  $(X, \mu_p)$, we mean that  $X$  is a smooth projective rational surface
with an effective  $\mu_p$-action.  
It is natural to assume that  $X$  is minimal in terms of  $G$-equivariant
birational morphisms (Definition 1.4).  

\par
We now state our results.
For  $X = {\bold F_e}$  in Theorem 1 (I), $X^{\mu_p}$  shuold be well-known
and is also determined in Lemma 4.3.

\par \vskip 1pc
{\bf Theorem 1.} {\it Let  $p$  be a prime number and let  $(X, \mu_p)$
be a minimal pair of a smooth projective rational surface and the group  $\mu_p$
acting effectively on  $X$.}

\par
(I) {\it If  $p$  is odd prime $($for  $p = 2$  see Theorem $4$ and Remark $5$ below$)$
and the  $\mu_p$-invariant sublattice  $(Pic X)^{\mu_p}$
has rank  $\ge 2$, then  $X$  is a Hirzebruch surface  $\bold F_e$ ($e \ne 1$)
and  $(X, \mu_p)$  is birationally equivariant to a pair  $(\bold P^2, \mu_p)$
given in Example $2.1$.}

\par
(II) {\it Suppose that  $(Pic X)^{\mu_p}$  has rank  $1$.
Then  $(X, \mu_p)$  is equal to one of the pairs in Examples $2.1-2.8$.
The fixed locus  $X^{\mu_p}$, $X$, $Y = X/\mu_p$,
Fano index  $r(Y)$, the types of all singularities on  $Y$
and the topological fundamental group of  $Y^0 = Y - \text{\rm Sing} \, Y$
are summarized in Table $1$ $($for odd prime  $p$  only$)$.}  

$$\text{\bf Table 1} \,\,\,\, 
\text{\rm [attached at the end of the paper]}$$

\par \vskip 1pc
{\bf Corollary 2.} {\it Let  $p$  be an odd prime number and let
$(X, \mu_p)$  be an arbitrary pair of a smooth projective
rational surface and the group  $\mu_p$  acting effectively on  $X$.
Then  $(X, \mu_p)$  is birationally equivariant to one of the pairs  
$(X_{\min}, \mu_p)$ in Table  $1$. Set $Y_{\min} = X_{\min}/\mu_p$.

\par
In particular, either  $X_{\min} = \bold P^2$, or
$Y_{\min} = \bold P^2$, or  $Y_{\min} = \overline{\bold F}_3$,
or  $(X_{\min}, \mu_p)$ $(p = 5)$  is one of the pairs  
in rows $5, 6, 7$ of Table $1$ given in Examples $2.4$ and $2.5$ 
$($see Lemmas $2.12-13$
for the uniqueness of the pairs in rows $6, 7)$.}

\par \vskip 1pc
In the result below, (2) is trivial, while
(1) is not so obvious; there is a rational surface
with at worst two (quotient) singularities such that its 
smooth part has infinite  $\pi_1$ (see [GZ3, \S 4]);
see also Remark 4.7.

\par \vskip 1pc
{\bf Corollary 3.} {\it Let  $(X, \mu_p)$  be as in Corollary $2$.
Set  $Y = X/\mu_p$  and  $Y^0 = Y - \text{\rm Sing} \, Y$.
Then  $\pi_1(Y^0)$  equals  $(1)$  or  $\mu_p$.
Moreover, when  $(X, \mu_p)$  is a minimal pair,
we have:}

\par
(1) {\it If  $X^{\mu_p}$  contains a curve, then  $Y^0$  is simply connected.}

\par
(2) {\it If  $X^{\mu_p}$  is a finite set, then  $\pi_1(Y^0)= \mu_p$
$($see also Lemma $4.4)$.}

\par \vskip 1pc
In the following, we denote by  $X^G = \{x \in X | g x = x$  for some  $1 \ne g \in G\}$
the fixed locus, $\sigma : X \rightarrow Y = X/G$  the quotient map, 
$Y^0 = Y - \text{\rm Sing} \, Y$  and  $r(Y)$  the Fano index.  

\par \vskip 1pc
{\bf Theorem 4.} {\it Let  $(X, G)$  be a minimal pair of a smooth projective
rational surface and an arbitrary finite group  $G$  
acting effectively on  $X$.
Then we have:}

\par
(I) {\it Suppose that the  $G$-invariant sublattice  $(Pic X)^G$
has rank $\ge 2$.  Then  $X$  has a  $G$-stable conic fibration
each singular fibre of which is a linear chain of two  $(-1)$-curves.}

\par
(II) {\it Suppose that  $(Pic \, X)^G$  has rank  $1$.
Then $X$  is a (smooth) del Pezzo surface and  $Y$  is a singular del Pezzo surface
with at worst quotient singularities so that  $\pi_1(Y^0)$  is finite;
one has  $\pi_1(Y^0) = G$  if the fixed locus  $X^G$  is a finite set.
Moreover, the following are true.}

\par
(1) {\it If  $r(Y) = 1$  and  $X^G$  is a finite set, then  $(X, G)$
is equal $($modulo $G$-equivariant isomorphism$)$ to one of the pairs in 
Examples $2.1b$ $(p=3)$, $2.5$, $2.9-11$.
$X^G, \, X, \, Y$ $= X/G, \, K_Y^2$  and the types of all singularities of  $Y$
are summarized in Table $2$ $($see Lemmas $2.13-15$ for the uniqueness of 
the pairs in the rows  $3, 4, 2, 6)$.}

\par
(2) {\it If  $r(Y) > 1$, then  $Y$  is either  $\bold P^2$  or the projective cone  
$\overline{\bold F}_e$ ($e \ge 2$)  with  $e | |G|$.}

$$\text{\bf Table 2} \,\,\,\, 
\text{\rm [attached at the end of the paper]}$$

\par \vskip 1pc
{\bf Remark 5.} (1) If  $G = \mu_2$  in Theorem 4 (I), then it is
birationally equivariant to some De Jonquieres involution of degree  $d \ge 2$
[BB];  when  $d = 2$, it is given in Example 2.6.

\par
(2) For a normal surface  $S$ (like  $Y$  in Theorem 4 (II)) 
with at worst quotient singularities,
$\bold Q$-ample anti-canonical divisor  $-K_S$  and  rank $Pic \, S = 1$,
the Fano index  $r(S) = 1$  holds if and only if  $S$  is a Gorenstein log del 
Pezzo surface other than the projective cone  $\overline{\bold F}_2$
(this cone has Fano index  $2$) (Lemma 1.9 and [MZ, Lemma 6]);
for such  $S$, it is also shown in [MZ, Lemma 6] that  $\pi_1(S^0)$ is 
abelian of order $\le 9$.

\par
(3) Kantor had classified automorphism groups of del Pezzo surfaces,
though it was not told which automorphism is lifted from a del Pezzo surface
of smaller degree.  In this sense, the result in Theorem 4 and Kantor's book 
together complete the picture of automorphism groups of rational surfaces.
In particular, we have to refer to [Kan] for the case where  $r(Y) < 1$.
See also  [MZ, Zh2].

\par
(4) The difference of our approach from others lie in two aspects:
(i) we determine also the fixed locus  $X^{\mu_p}$  and the quotient surface
$X/\mu_p$ and (ii) we include both the geometric approach (bottom up), and
the algebraic approach as an Appendix, though the uniqueness
and realizability of pairs are only treated in the geometric approach.

\par \vskip 1pc
{\bf Acknowledgement.}  This work was done during the author's visit at
University of Michigan in 1999 who would like to thank its Department of
Mathematics for the hospitality. 
We would also like to thank Jonghae Keum for very valuable discussion.
This work was partially financed by an Academic Research Grant
of National University of Singapore.

\par \vskip 2pc
{\bf \S 1. Preliminary results}

\par \vskip 1pc
{\bf 1.1.} Let  $(X, G)$  be a pair of a smooth rational projective surface
and a non-trivial finte group  $G$  acting effectively on  $X$.
Denote by  $Y = X/G$  the quotient 
surface and  $\sigma : X \rightarrow Y$  the quotient map.
Let  $f : Z \rightarrow Y$  be the minimal resolution.
Note that $Y$ is a rational surface by Luroth's theorem
and $Y$ has at worst quotient singularities
and hence is simply connected [Ko, Th 7.8].

\par
We assert that $X^G$ is non-empty. Indeed, if $X^G$
is empty then the quotient map $X \rightarrow Y$
is an unramified finite morphism of degree $|G|$ over the simply-connected
surface $Y$, whence $|G| = 1$, a contradiction.

\par
The following is well known [Bri, Satz 2.11]; 
for the smoothness of  $X^{\mu_p}$  we diagonalize the action
locally and see that the  $\mu_p$  fixed part is defined by
a local coordinate (the eigenvector w.r.t. to the eigenvalue $\ne 1$).

\par \vskip 1pc
{\bf Lemma 1.2.} (1) {\it The fixed locus  $X^G$  is non-empty.
If  $G = \mu_p$  then  $X^{\mu_p}$  is a disjoint union of
smooth curves  $R_i$  and finitely many points  
$p_j$ ($1 \le j \le s$; $s \ge 0$).}

(2) {\it The surface  $Y = X/G$  is a  ${\bold Q}$-Gorenstein
normal rational surface with singularities.
If  $G = \mu_p$, then  $q_i := \sigma(p_i)$  is
a cyclic quotient singularity of type  $\frac{1}{p}(1, k_i)$  
for some  $1 \le k_i \le p-1$;
one has  $\text{\rm Sing} \, Y = \{q_1, \cdots, q_s\}$  and  $\sigma^{-1}(q_i) = p_i$.}

\par
(3) {\it Suppose that  $G = \mu_p$.  Then
$Y = X/\mu_p$  is Du Val at  $q_i$ (i.e., Gorenstein in the present quotient singularity case) 
if and only if  $k_i = p-1$  (this is always true when  $p = 2$).
In general, $p K_Y$  is Cartier.}

\par
(4) {\it The quotient map  $\sigma$  is unramified outside the fixed locus
$X^G$.  If  $G = \mu_p$, the ramification formula has the form
($\bold Q$-linear equivalence) :
$K_X \sim_{\bold Q} \sigma^*(K_Y) + (p-1) \sum_i R_i$.}

\par
(5) {\it The  $\sigma$-invariant sublattice  $(Pic \, X)^G \otimes {\bold Q}$
has rank equal to the Picard number  $\rho(Y) =$ rank $Pic \, Y$  
of the quotient surface $Y$.}
 
\par \vskip 1pc 
{\bf 1.3.} On surfaces, quotient singularity and log terminal singularity
are equivalent [Kaw, Cor 1.9].  
So there is a  ${\bold Q}$-effective divisor  $\Delta$
supported on  $f^{-1}(\text{\rm Sing} \, Y)$  and with the integral part  $[\Delta] = 0$,
such that
$$K_Z = f^*(K_Y) - \Delta.$$
Write  $\Delta = \sum_i \Delta_i$  where  $\Delta_i$  is supported
on  $f^{-1}(q_i)$.  Then  $\Delta_i = 0$  if and only if
$q_i$  is Du Val.

\par \vskip 1pc
{\bf 1.4 Definition.} 
Fix a group  $G$.  Let  $(X_i, G)$  be two pairs where
$G$  acts effectively on   $X_i$.  A {\it $G$-equivariant birational morphism}  
$\tau : (X_1, G) \rightarrow (X_2, G)$  
is a birational morphism  $\tau : X_1 \rightarrow X_2$,
satisfying  $\tau(gx) = g \tau(x)$  for every  $g \in G$.
The existence of such  $\tau$  is equivalent to that of a  $G$-stable
divisor on  $X_1$  which can be smoothly blown down.
If the  $G$-invariant sublattice  $(Pic \, X)^G \otimes \bold Q$
has rank 1, then there is no such  $\tau$  and  $(X, \mu_p)$
is minimal in the sense below.

\par
Two pairs  $(X_i, G)$  are birationally  $G$-equivariant
if there is a birational map  $X_1 \cdots \rightarrow X_2$
which can be decomposed as  $f_1 \circ \cdots \circ f_n$
such that for each  $i$  either  $f_i$  or  $f_i^{-1}$
is a  $G$-equivariant biraitonal morphism.

\par
$(X, G)$  is called a {\it minimal pair},
if for any  $G$-equivariant birational morphism  
$\tau : (X, G) \rightarrow (X_2, G)$, one has  $\tau = $id.

\par
Let  $(X, G)$  be a pair (with  $X$  rational and  $G$  finite) and  let  $Y = X/G$.
Suppose that  $-K_Y$  is  {\bf Q}-ample.  Wrtie  $-K_Y \sim_{\bold Q} r P$,
where  $r$  is a positive rational number and  $P$  a Cartier ample divisor.
Let  $r(Y)$  be the largest (hence  $P$  is the
``smallest'') among such expression, noting that
$Pic \, Y$  is a torsion free  $\bold Z$-module of finite rank
($Y$  is simply connected).
By the same reasonning, the divisor class of
$P$  is uniquely determined by  $-K_X$  or  $X$.
This  $r(Y)$  is called the {\it Fano} index of  $Y$.
When  $G = \mu_p$, one can write
$r(Y) = m/p$  with a positive integer  $m$  because  $p K_Y$  is Cartier.

\par \vskip 1pc
{\bf Remark 1.5.} 
If  $X$  is a smooth Fano $n$-fold (i.e., $-K_X$  is ample)
then  $r(X) \le n + 1$, and  $r(X) = n$
(resp. $r(X) = n+1$)  if and only if  $X$  is a smooth quadric hypersurface
in  $\bold P^{n+1}$  (resp. $X = {\bold P}^n$) [KO].

\par
Let  $\overline{\bold F}_e$ ($e \ge 2$)  be the projective cone,
with vertex  $q_1$, over a (smooth) rational curve of degree  $e$  in  $\bold P^e$.
Then the resolution of the vertex is the Hirzebruch surface
$\bold F_e$ (see [Hart]), where the  $(-e)$-curve is
the inverse of the vertex.  $\overline{\bold F}_e$  can also be
embedded into  $\bold P^{e+1}$  as a non-degenerate
surface of degree  $e$ (see [Nag]).
The hyperplane section  $\overline H$  
of  $\overline{\bold F}_e \subseteq \bold P^{e+1}$  
is the generator of the Picard lattice
and is the image of a section on  $\bold F_e$  disjoint from the  $(-e)$-section;
so  $\overline H^2 = e$.  One sees that  $r = (e+2)/e > 1$.

\par \vskip 1pc
{\bf 1.6.} Suppose that  $G = \mu_p$. 
The induced  $\mu_p$  action on  $Pic X \otimes \bold C$  can be 
diagonalized.  Since  $\mu_p$  acts on the integral lattice  $Pic \, X$,
there is generator  $h$  of  $\mu_p$  satisfying, 
where  $\zeta_p = exp(2 \pi \sqrt{-1}/p)$ :
 
$$h^*|Pic X \otimes \bold C = \text{\rm diag} [I_c, M_p^{\oplus t}], \,\,\,\,
M_p = [\zeta_p, \zeta_p^2, \dots, \zeta_p^{p-1}].$$

\par \vskip 0.5pc
{\bf Lemma.} (1) {\it $\text{\rm rank (Pic} X)^{\mu_p} \otimes \bold Q = c$;
$K_X^2 = 10 - c - (p-1)(2+c - s - \sum_i (2 - 2g(R_i))$.}

\par
(2) {\it Writing  $X^{\mu_p} = \coprod R_i \coprod \{p_1, \dots, p_s\}$,
we have  $s + \sum_i (2 - 2g(R_i)) = c + 2 - t$.}

\par
(3) {\it Let  $k_i$ $(1 \le i \le p-1)$ be the number of 
isolated  $\mu_p$-fixed points at which a generator of
$\mu_p$  can be diagonalized as  $(\zeta, \zeta^i)$
with  $\zeta$  a primitive  $p$-th root of  $1$.
Then  
$$\gather 1 = \sum_j (\frac{1-g(R_j)}{2} + \frac{(p+1)R_j^2}{12}) +
\frac{1}{p-1} \sum_{i=1}^{p-1} \sum_{\zeta} k_i/(1-\zeta)(1-\zeta^i) = \\
\sum_j (\frac{1-g(R_j)}{2} + \frac{(p+1)R_j^2}{12}) +
k_1\frac{5-p}{12} + k_2\frac{11-p}{24} + k_3 a_3 + k_4 a_4 
+ \cdots \endgather $$
where  $\zeta$  runs over the set of primitive  $p$-th root of $1$,
where  $a_3 = 1/4, \, a_4 = 1/2$  when  $p = 5$.}

\par \vskip 1pc
{\it Proof.} Applying the topological fixed-point formula,
we obtain  
$$s + \sum_i (2 - 2g(R_i)) = \chi_{top}(X^{h}) = \sum_i (-1)^i Tr(h^*|H^i(X, {\bold C}))
= 2 + c + t \, Tr(M_p) = 2 + c - t.$$
The Picard number  $\rho(X)$  is  $10 - K_X^2$  and also equals  $c + t(p-1)$.
So (1) and (2) are proved.

\par
By the holomorphic Lefschetz fixed point formula [ASIII, p.567],
one has
$$\gather
1 = \sum_i (-1)^i Tr (h^*|H^i(X, \Cal O_X)) = \\
\sum_i 1 / \det(1-h|T_{p_i})
+ \sum_j (1- g(R_j))/(1 - \zeta_p^{n_j}) - 
\sum_j R_j^2 \zeta_p^{n_j}/(1 - \zeta_p^{n_j})^2,
\endgather $$
where  $T_{p_i}$  is the tangent space of  $X$  at  $p_i$,
$h$  is the generator of  $\mu_p$  and  $h^*$  acts on
the normal bundle of  $R_j$  by a multiple  $\zeta_p^{n_j}$ 
(a primitive $p$-th root of 1).
Letting  $h$  run in the set of generators of  $\mu_p$  and taking
sums for both sides of the above equality, to prove (3) we only need to show
$$(p-1)/2 = \sum_{i=1}^{p-1} 1/(1-\zeta_p^i), \,\,\,\,
(1-p^2)/12 = \sum_{i=1}^{p-1} \zeta_p^i/(1-\zeta_p^i)^2.$$
These can be checked by using  $p = (1-x)(x^{p-2} + 2 x^{p-3} +
\cdots + (p-2) x + p-1)$  with  $x = \zeta_p$
to get rid of the denominators.  The equalities above were 
originally calculated by Cay Horstman and was kindly brought to our attention
by Jonghae Keum.

\par \vskip 1pc
In Lemmas 1.7-1.10 below, except Lemma 1.7 (4) and Lemma 1.8, we assume only
that  $X$  is a smooth rational surface and  $G$  a non-trivial finite group
acting on it.

\par \vskip 1pc
{\bf Lemma 1.7.} {\it Suppose that  
$\text{\rm rank} (Pic \, X)^G \otimes {\bold Q} = 1$.  
Then we have:}

\par
(1) {\it $X$  is a del Pezzo surface.
Hence  $d = K_X^2$  satisfies
$1 \le d \le 9$.  $d = 9$  if and only if  $X = {\bold P}^2$;
$d = 8$  if and only if  $X$  is the Hirzebruch surface  ${\bold F}_e$
with  $e = 0$  or  $e = 1$; $d \ge 2$  if and only if  $X$  is the blow-up
of  ${\bold P}^2$  at  $9-d$  points in general position.}

\par
(2) One has  $Pic \, Y = \bold Z P$ (see Definition $1.4$).
If  $K_X^2 \le 7$, then  $(Pic \, X)^G = \bold Z K_X$.

\par
(3) {\it $-K_Y$  is  {\bf Q}-ample.  A general member of
$|P|$  is smooth and irreducible, which does not 
pass through the singular locus of  $Y$; one has also
$2g(P) - 2 = P . (K_Y + P)$.}

\par
(4) {\it Suppose further that  $G = \mu_p$.
Then  $X^{\mu_p}$  is either a finite set, or a union of a
smooth irreducible curve  $R$  and finitely many points.}

\par \vskip 1pc
{\it Proof.} Clearly, both the pull-back  $H$  on  $X$  of an ample divisor on  $Y$
and  $-K_X$  are generators of the rank one  ${\bold Q}$-module
$(Pic X)^G \otimes {\bold Q}$.  Noting that the Kodaira dimension
$-\infty = \kappa(X) < 2$, $-K_X$  is a positive multiple of  $H$
and (1) follows [Man2].

\par
The first part of (2) is true because  $Pic \, Y$  is a rank
one lattice and $-K_Y$ is ample (see (3)).  
Let  $C$  be any  $G$-stable Cartier divisor.
Then  $C = (m/n)(-K_X)$  for some coprime positive integers.
Intersecting this with a  $(-1)$-curve  $E$  on  $X$, one obtains
$n (C . E) = m$  and  $n | m$, whence  $n = 1$  and (2) is proved.

\par
For (4), if  $X^{\mu_p}$  contains two (disjoint) curves  $R_1, R_2$
then both  $R_i$  are positive multiples of  $H$  and this leads to 
that  $0 = R_1 . R_2 = H^2 \times$ (a positive number), a contradiction.

\par
For (3), since  $\rho(Y) = 1$, either  $K_Y$  or  $-K_Y$  is  {\bf Q}-ample.
By the ramification formula (similar to the one in Lemma 1.2) and the fact that the Kodaira
dimension  $\kappa(X, K_X) = -\infty$, we see that  $\kappa(Y, K_Y) = - \infty$;
so  $-K_Y$  is  {\bf Q}-ample.

\par
By the main theorem of [Am] or [Alex], dim $|P| \ge 1$
and (a general member)  $P$  is smooth irreducible.  Since  $P$  is Cartier
and  $Y$  has at worst rational singularities (for the second equality) 
and by the Riemann-Roch theorem (for the third) one has (cf. [Art, Th 2.3]):
$$\gather 
\chi({\Cal O}_P) = \chi({\Cal O}_Y) - \chi({\Cal O}_Y(-P)) =
\chi({\Cal O}_Z) - \chi({\Cal O}_Z(-f^*P)) =
\frac{-1}{2}f^*P . (K_Z + f^*P) \\ 
= \frac{-1}{2}(P' + D) (K_Z + P') = \chi({\Cal O}_{P'}) + \frac{-1}{2} D . (K_Z + P'),
\endgather $$
where  $P'$  is the proper transform on  $Z$  of  $P$  and  $D$  an 
effective divisor with support in  $f^{-1}(\text{\rm Sing} \, Y)$  
such that  $f^*P = P' + D$.
Since  $P \cong P'$, we get  $0 = D . (K_Z + P') \ge D . P'$,
where the inequality is true because each component of  $D$
is a  $(-n)$-curve with $n \ge 2$ [Bri, Satz 2.11].  Thus  $D^2 = D . f^*P
= 0$  and  $D = 0$  for  $D$  is negative definite [Mum, page 230].
So  $f^*P = P'$  and  $P$  is away from the singular locus of  $Y$.
Now  $2g(P) - 2 = 2g(P') - 2 = f^*P . (K_Z + f^* P) = P . (K_Y + P)$
because  $f_*K_Z = K_Y$.

\par \vskip 1pc
{\bf Lemma 1.8.} {\it Suppose that the quotient surface  $Y = X/\mu_p$  satisfies
$\rho(Y) = 1$  and   $r(Y) < 1$.
Then the fixed locus  $X^{\mu_p}$  is a finite set.}
\par \vskip 1pc
{\it Proof.} Write  $-K_Y = r P$  with  $r = r(Y) < 1$.  
Suppose the contrary that  $X^{\mu}$  contains an irreducible curve  $R$.
Then  $X^{\mu}$  is a union of  $R$  and points  $p_i$ (Lemmas 1.2 and 1.7).
Note that  $B = \sigma(R)$  is away from  $\text{\rm Sing} \, Y$ (Lemma 1.2),
is Cartier and satisfies  $\sigma^* B = p R$.
Write  $B = b P$, where  $b \ge 1$  by the maximality of  $r$.
Then the ramification formula implies that
$$-K_X = -(\sigma^*K_Y + (p-1)R) = [rp - b(p-1)]/p \sigma^*P.$$
Since  $-K_X$  is ample (Lemma 1.7),  $r > b(p-1)/p \ge (p-1)/p$.
This and the fact that  $r = m/p$  with an integer  $m$
would imply that  $m \ge p$  and  $r \ge 1$.  This contradiction
proves the lemma.

\par \vskip 1pc
The following two results are essentially proved in [Fuj, Chapter 1, \S 5].
For the convenience of
the reader we give a kind of new proof here.

\par \vskip 1pc
{\bf Lemma 1.9.} {\it Suppose that the quotient surface  $Y = X/G$  satisfies
$r(Y) = 1$.  Then  $Y$  is Du Val (i.e., Gorenstein in the present case) 
and a general member of  $|-K_Y|$  is a smooth
elliptic curve which does not pass through the singular locus of  $Y$.}

\par \vskip 1pc
{\it Proof.} In view of Lemma 1.7, we only need to
show that  $Y$  is Gorenstein and  $P \sim -K_Y$.
By the assumption  $K_Y + P$  is  $\bold Q$-linearly
equivalent to zero.  Then  $0 \sim_{\bold Q} f^*(K_Y + P) = (K_Z + P') + \Delta$;
here  $P' = f^*P$  is a smooth curve with
$p_a(P') = 1 > 0$ (cf. 1.3 and 1.7), and hence the Riemann-Roch
theorem implies that  $|K_Z + P'| \ne \emptyset$.  Thus
$\Delta = 0$, whence  $Y$  has only Du Val singularities (cf. 1.3).
Finally, since the two Cartier divisors  $-K_Y$  and  $P$  are  
$\bold Q$-linear equivalent, they are linear equivalent
because the rational surface  $Y$  is simply connected
and hence  $Pic \, Y$  is torsion free.  The lemma is proved.

\par \vskip 1pc
{\bf Lemma 1.10.} {\it Suppose that the quotient surface  $Y = X/G$  satisfies
$r = r(Y) > 1$.  Then (a general member)  $P$  is a smooth rational curve
away from the singular locus of  $Y$.  Moreover, $(r-1)P^2 = 2$.}

\par \vskip 1pc
{\it Proof.} Substituting  $-K_Y = rP$  into the equality
in Lemma 1.7, we get  $2g(P) - 2 = (1-r)P^2 < 0$.  Thus  
$g(P) = 0$  and the current lemma follows from Lemma 1.7.

\par \vskip 1pc
{\bf \S 2. Examples}

\par \vskip 1pc
In this section, we shall construct examples of pairs  $(X, \mu_p)$
(see Theorems 1 and 4 in the Introduction).

\par \vskip 1pc
{\bf 2.1.} Suppose that  $\mu_p$  acts effectively on  $X = {\bold P^2}$
with homogeneous coordinates  $X, Y, Z$.
Then one can diagonalize a suitable generator  $g$  of  $\mu_p$  as one of the following,
where  $\zeta_p = exp(2 \pi \sqrt{-1}/p)$:
$$\text{\bf 2.1a:} \,\,\,\, g = \text{\rm diag} [1, 1, \zeta_p]; \,\,\,\,\,
\text{\bf 2.1b:} \,\,\,\, \text{\rm diag} [1, \zeta_p, \zeta_p^v] \,\, (2 \le v \le p-1).$$

\par
In 2.1a, $X^{\mu}$  is a union of the line  $Z = 0$
and the point  $p_1 = [0, 0, 1]$.  This  $p_1$  dominates a singularity
$q_1$  of   $Y := X/\mu_p$  of type  $\frac{1}{p}(1,1)$.
It is easy to see that  $Y$  is the projective cone  $\overline{\bold F}_p$
with the vertex at  $q_1$.  The  $\bold Z/(p)$-covering map  
$\sigma : X \rightarrow Y$  is branched along the vertex
and a smooth hyperplane  $B$ ($\sim \overline H$  in notation of Remark $1.5$).
One has  $r(Y) = (p+2)/p$ (Remark 1.5).

\par
In 2.1b, one must have  $p \ge 3$  and  $X^{\mu}$  is a union 
of three points  $p_1 = [1, 0, 0], \, p_2 = [0, 1, 0], \, p_3 = [0, 0, 1]$.
These  $p_i$  dominate singular points  $q_i$
of  $Y := {\bold P}^2/\mu_p$.  The  $q_i$  are respectively of type  
$\frac{1}{p}(1, v)$, $\frac{1}{p}(1, p+1-v)$  and  $\frac{1}{p}(1, u)$  
with  $uv \equiv v-1$ (mod $p$).   
One has  $\pi_1(Y^0) = \mu_p$ (Lemma 4.4).
One sees also that  $Y$ is Du Val if and only if  $p = 3$; if this is the case
then  $g = \text{\rm diag} [1, \zeta_3, \zeta_3^2]$  and
$Y$  is a Gorenstein log del Pezzo surface of rank
1 with 3 type  $\frac{1}{3}(1,2)$  singularities
and also  $r(Y) = 1$ [MZ1, Lemma 6].  If  $p \ge 5$, then  $r(Y) = 1/p, \, 3/p$ 
(Lemma 1.9 and Propositions 3.1 and 3.3).

\par \vskip 1pc
{\bf 2.2.} Let  $Y = \bold P^2$, $p = 3$  and  $\sigma : X \rightarrow Y$  
the triple cover totally branched along a smooth plane cubic  $B$; 
$X$  is a del Pezzo surface
with  $K_X^2 = 3$  and  $R \in |-K_X|$, where  $\sigma^*B = 3 R$.

\par \vskip 1pc
{\bf 2.3.} Let  $Y = \overline{\bold F}_p$, $p = 3$  and  $\sigma : X \rightarrow Y$ 
the the triple cover branched along a smooth genus-$2$  curve  $B$
and the vertex  $q_1$ ($\notin B$) of the cone  $Y$.  Then  $X$  is a del Pezzo surface
with  $K_X^2 = 1$  and  $|-K_X|$  contains 6 members of cuspidal 
rational curves lying over the 6 generating lines of the cone  $Y$
tangent to the branch curve  $B$.  Blowing up the unique
base point of  $|-K_X|$ [Dem, Proposition 2, page 40] with  $E$  the exceptional curve,
one gets a relatively minimal rational elliptic surface  
$\varphi : \widetilde X \rightarrow \bold P^1$  with a section  $E$
and six type  $II$  singular fibres; so the Mordell Weil group
of the fibration   $\varphi$  is torsion free and of full rank 8.
There is an induced  $\mu_3$-action on  $\widetilde X$  fixing (point wise)
the section  $E$.

\par \vskip 1pc
{\bf 2.4} (the rows 5 and 6 of Table 1). 
Here we construct a 1-dimensional family  $(X_s, \mu_5)$
($s \in \bold P^1 \setminus \{0, \pm 1\}$)
and a unique pair  $(X_{II}, \mu_5)$, where each surface is
a degree 1 del Pezzo surface on which  $\mu_5$  acts effectively
and fixes (point wise) a smooth member in the anti-canonical 
linear system.  When  $X = X_s$ (resp. $X = X_{II}$),
$|-K_X|$  has  $10$  nodal members forming two  $\mu_5$-orbits, 
and one cuspidal member (resp. $5+1$  cuspidal members forming
two  $\mu_5$-orbits).

\par
Let  $\widetilde Z \rightarrow \bold P^1$  with  $\widetilde Z = \widetilde Z_I$
(resp. $\widetilde Z = \widetilde Z_{II}$) be the unique
ellipitc surface with (only one) section  $E$, a type  $II^*$  fibre 
$\widetilde Z_{t=0}$, and two type  $I_1$  fibres at  $t = \pm 1$
(resp. a single type  $II$  fibre at  $t = \infty$) [MP, Th 5.4].
Express the type  $II^*$  fibre as  $D_1 + 5C + D_2$, where  $C$  is 
a $(-2)$-curve, and  Supp $D_i$  are {\it the} two disjoint chains of  $(-2)$-curves
of length 4 so that the section  $E$  meets a tip component of  $D_2$.
Let  $\widetilde Z \rightarrow Y$ (late on, referred as  $Y = Y_I, Y = Y_{II}$
respectively) be the contraction of  $E + D_2, \, D_1$  to a smooth point  
$q$  and a type  $\frac{1}{5}(1, 4)$  
singular point  $q_1$; both points lie on the image of  $C$ (also denoted by  $C$).  

\par
Let  $B$  be the image on  $Y$  of a smooth fibre
$\widetilde Z_{t=s}, s \ne 0, \pm 1$ (resp. $\widetilde Z_{t=1}$)
when  $Y = Y_I$ (resp. $Y = Y_{II}$).
Since fibres on  $\widetilde Z$  are linearly equivalent, pushing down,
we get an induced relation  ${\Cal O}_Y(C)^{\otimes 5} \cong {\Cal O}_Y(B)$.
This gives rise to a  $\mu_5 \cong \bold Z/(5)$-Galois cover :
$$\sigma : X = Spec \oplus_{i=0}^4 
{\Cal O}_Y(-i \, C) \longrightarrow Y,$$
(referred as  $X = X_s, X = X_{II}$  respectively)
which is etale outside the smooth elliptic  $B$  in  $|-K_Y|$  
and the only singularity  $q_1$  of  $Y$; along  $B$
the map  $\sigma$  is totally branched.

\par \vskip 1pc
{\bf 2.4.1.} Conversely, suppose that  $5 \overline{M} \sim B$  is a relation on
$Y$  with  $\overline{M}$  a Weil divisor.  We now show that
${\Cal O}_Y(\overline{M}) \cong {\Cal O}_Y(C)$.
Pulling back it by the
minimal resolution  $f : Z \rightarrow Y$ ($\widetilde Z \rightarrow Z$
is the contraction of  $E + D_2$), one has  $5M + D_1' \sim B$  on  $Z$,
where  $M$  is mapped to  $\overline{M}$, the same  $B$  denotes
its preimage on  $Z$  and  $D_1'$  is supported on the support
of (the image on  $Z$  of)  $D_1$.  On  $Z$, (the image of) $C$
satisfies  $C . B = 1$, while each component of  $D_1$  has zero
intersection with  $B$.  Using these to intersect
the above relation with  $C$  and components of  $D_1$,
one sees that  $D_1' - D_1 = 5D$  with  $D$  a  $\bold Z$-combination of
irreducible components of  $D_1$.  Thus  $5(M + D) + D_1 \sim B$  on  $Z$.
On the other hand, one has  $5C + D_1 \sim B$  on  $Z$.
These relations imply that  $5(M + D - C) \sim 0$  and hence
$M + D \sim C$  for the rational surface  $Z$  has torsion free
$Pic \, Z$.  Passing to  $Y$, one gets  ${\Cal O}_Y(\overline{M}) \cong {\Cal O}_Y(C)$. 

\par \vskip 1pc
It is easy to see that  $X$  is smooth and
$-K_X = \sigma^*(-K_Y - 4C) = \sigma^*(5C - 4C) = \sigma^*C$
so that  $-K_X$  is nef and big with  $K_X^2 = 1$  
because  $C^2 = 1/5$  on  $Y$;
every member  $F$ ($\ne B, \, 5C$)  in  $|-K_Y|$
has total transform on  $X$  splitting into
5 elliptics meeting at the unique point lying over  $q$ ($= B \cap C$),
while  $\sigma^*B = 5R$  with a smooth elliptic  $R \in |-K_X|$
and  $\sigma^*C$  is a cuspidal curve in  $|-K_X|$  with a cusp at the
point  $p_1 = \sigma^{-1}(q_1)$;
thus  $\sigma$  is totally ramified exactly  along  $R$  and the point  $p_1$,
and  $X$  is a del Pezzo surface with  $X^{\mu_5} = R \coprod \{p_1\}$;
indeed, $X$  has no  $(-2)$-curves and the only singular members
($\ne \sigma^*C$)  in  $|-K_X|$  are $10$ nodal curves lying over
the two type  $I_1$  fibres on  $\widetilde Z_I$
(resp. five cuspidal curves lying over the type  $II$  fibre on  $\widetilde Z_{II}$).

\par
Finally, we have  $\text{\rm rank} \, Pic \, Y = 1$  and  $K_Y^2 = 5$
(noting that  $Y$  is Du Val).  One has also  $r(Y) = 1$
and  $\pi_1(Y - \text{\rm Sing} \, Y) = (1)$ [MZ, Lemma 6].

\par \vskip 1pc
{\bf 2.5} (the row 7 of Table 1).  We shall calculate  $X^{\mu_5}$ 
and determine the type of singularities of  $Y = X/\mu_5$
for the unique pair  $(X, \mu_5)$, where $X$ is the unique 
del Pezzo surface with  $K_X^2 = 5$ (Lemma 2.13).

\par
Let  $\widetilde X \rightarrow \bold P^1$  be a relatively
minimal elliptic surface with a section and with two type $I_5$
fibres  $F_1, F_2$; such  $\widetilde X$  is unique and the fibration has
the Mordell Weil group  $\bold Z/(5)$ [MP, Th 5.4]; 
it has also two type  $I_1$
singular fibres.  Using Shioda's height pairing
[Sh, Th 8.6], one can verify that the five sections  $P_i$ ($0 \le i \le 4$)
are disjoint and meet distinct fibre components of  $F_1, F_2$.
Let  $\widetilde X \rightarrow X$  be the blow-down of the 
sections  $P_i$.  Then  $X$  is the del Pezzo surface with  $K_X^2 = 5$.
The translation automorphism given by the section  $P_1$
induces an automorphism  $g$  on  $X$  so that  $X^{\mu_5}$  consists of
two points (the images of the nodes of the two type  $I_1$
fibres on  $\widetilde X$), where  $\mu_5 = \langle g \rangle$.
Set  $Y = X/\mu_5$  and  $Y^0 = Y - \text{\rm Sing} \, Y$.
Then  $Y$  is a Gorenstein log del Pezzo surface with two type
$\frac{1}{5}(1, 4)$  singularities, $K_Y^2 = 1$, rank $Pic \, Y = 1$,
and  $r(Y) = 1$ (see [MZ, Lemma 6]).

\par
Conversely, one can show that such kind of singular del Pezzo surface 
$Y$  is unique modulo isomorphism and that  $\pi_1(Y^0) = \mu_5$;
hence the quotient map  $\sigma : X \rightarrow X/\mu_5$  here,
is the completion for the universal covering map of 
the smooth part of the {\it unique} Gorenstein log del Pezzo 
surface with two type  $\frac{1}{5}(1, 4)$  singularities
(see Lemma 4.4).

\par \vskip 1pc
{\bf 2.6.} $\mu_p$ ($p = 2$) acts on  $X = \bold P^1 \times \bold P^1$
by  $(x, y) \mapsto (y, x)$ (switching the fibrations).
One sees that  $Y = X/\mu_p$  is  $\bold P^2$
and the quotient map  $\sigma : X \rightarrow Y$  is
branched along a smooth conic, whose inverse on  $X$  is the diagonal. 
This  $\mu_p$  is birationally equivalent to {\it De Jonquieres involution} of 
degree 2 [BB, Example 1.6].

\par \vskip 1pc
{\bf 2.7.} Let  $Y = \bold P^2$, $p = 2$  and  $\sigma : X \rightarrow Y$  
the double cover branched along a smooth quartic curve.  $\mu_p$
is called {\it Geiser's involution} on the del Pezzo surface  $X$
with  $K_X^2 = 2$.  Conversely, if  $X$  is a del Pezzo surface with
$K_X^2 = 2$  then  $\Phi_{|-K_X|}$  is the  $\sigma$  above [Dem, Chapter V, \S 4].

\par \vskip 1pc
{\bf 2.8.} Let  $Y = \overline{\bold F}_p$ ($p = 2$)  and  
$\sigma : X \rightarrow Y$  
the double cover branched along a smooth genus-$4$ curve  $B$  and the 
vertex  $q_1$ ($\notin B$) of the cone  $Y$.   
Then  $X$  is a del Pezzo surface
with  $K_X^2 = 1$.  $\mu_p$  is called the {\it Bertini involution.}
Conversely, if  $X$  is a del Pezzo surface with
$K_X^2 = 1$  then  $\Phi_{|-2K_X|}$  is the  $\sigma$  above [Dem, Chapter V, \S 5].

\par \vskip 1pc
{\bf 2.9} (the row 4 of Table 2).
We construct a pair  $(X, \mu_6)$  with  
$X$  the (unique) del Pezzo surface with  $K_X^2 = 6$  and the group  
$\mu_6$  acting effectively on  $X$  such that  $X^{\mu_6}$  is 
a finite set and  $Y = X/\mu_6$  has exactly  $3$  singularities 
of type  $\frac{1}{2}(1,1), \, \frac{1}{3}(1,2), \,
\frac{1}{6}(1,5)$ (types $A_1, \, A_2, \, A_5$ in other notation) 
as all of its singularities.

\par
Let  $\widetilde X \rightarrow \bold P^1$  be a relatively 
minimal rational elliptic fibration with type  $I_1, I_2, \, I_3, \, I_6$
singular fibres  $F_0, F_1, F_2, F_3$  and a section  $P_0$.  
Such an elliptic surface is unique
[MP, Th 5.4].  Write  $F_1 = C_0 + C_1, \, F_2 = \sum_{i=0}^2 D_i, \,
F_3 = \sum_{j=0}^5 E_j$  with  $D_i . D_{i+1} = E_j . E_{j+1} = 1$
so that  $P_0$  meets  $C_0, \, D_0, \, E_0$.
Using the height-pairing in [Sh, Th 8.6], one sees that the Mordell Weil
group of the fibration is  $\bold Z/(6)$  and its generator
$P_1$  meets  $C_1, \, D_1, \, E_1$  after suitable relabelling.

\par
Denote by  $g$  the translation automorphism given by the section
$P_1$.  One sees that  ${\widetilde X}^{\mu_6}$, with  
$\mu_6 = \langle g \rangle$, consists of 6 points : $g$ (resp. $g^3$; $g^2$) 
fixes the node  $p_1$ (resp. the two nodes  $p_2, \, p_3$; 
the three nodes  $p_4, \, p_5, \, p_6$) of the type  $I_1$ 
(resp. $I_2$; $I_3$) fibre.  Let  $\widetilde X \rightarrow X$  
be the smooth blow-down of the six (disjoint) sections so that  
$X$  is the degree 6 del Pezzo surface.  Then  $X$  has an induced
$\mu_6$  actions so that  $X^{\mu_6}$  is still a 6-point set
(the image of  ${\widetilde X}^{\mu_6}$).  One sees that
$Y = X/\mu_6$  has exactly three singular points of type  $A_5$, $A_2$, $A_1$
(the images of  $p_1, \{p_2, p_3\}, \{p_4, p_5, p_6\}$),
rank $Pic \, Y = 1$  and  $K_Y^2 = 1$.

\par \vskip 1pc
{\bf 2.10} (the row 5 of Table 2). 
Suppose that  $G \cong (\bold Z/(3))^{\oplus 2}$
acts effectively on  $\bold P^2$  with coordinates  $X, Y, Z$,
so that  $X^G$  is a finite set.  Take generators  $g_1, \, g_2$  of  $G$.
By the assumption on  $X^G$, one has
$g_1 = \text{\rm diag} [1, \zeta_3, \zeta_3^2]$  after
a change of coordinates.  Now the commutativity of  $g_1, \, g_2$  
in  $PGL_2(\bold C)$  implies that  $g_2 = (a_{ij})$  with
$a_{13} a_{21} a_{32} = 1$  and  $a_{ij} = 0$  for other entries.
One sees easily that each order 3 subgroup of  $G$  fixes exactly
3 points, and  $|X^G| = 12$.  For instance, $X^{g_1} =
\{[1, 0, 0], \, [0, 1, 0], \, [0, 0, 1]\}$
and  $X^{g_2} = \{[1, a_{21} \zeta_3^i, \zeta_3^{2i}/a_{13}] |
0 \le i \le 2 \}$.

\par \vskip 1pc
{\bf 2.11} (the rows 2 and 6 of Table 2).  We shall construct: 

\par
(1) a pair  $(X, \mu_4)$  with  
$X = \bold P^1 \times \bold P^1$
and the group  $\mu_4$  acting effectively on  $X$
such that  $X^{\mu_4}$  is a finite set and  $Y = X/\mu_4$  has exactly
$3$  singularities of type  $\frac{1}{2}(1,1), \, \frac{1}{4}(1,3), \,
\frac{1}{4}(1,3)$ (i.e., types $A_1, A_3, A_3$) as all of its singularities; and

\par
(2) a pair  $(X, G)$  with  
$X = \bold P^1 \times \bold P^1$
and the group  $G \cong \mu_2 \times \mu_4$  acting effectively on  $X$
such that  $X^G$  is a finite set and  $Y = X/G$  has exactly
$4$  singularities of type  
$\frac{1}{2}(1,1), \, \frac{1}{2}(1,1), \, \frac{1}{4}(1,3), \,
\frac{1}{4}(1,3)$ as all of its singularities.

\par
Let  $\widetilde X \rightarrow \bold P^1$  be a relatively 
minimal rational elliptic fibration with type  $I_2, \, I_2, \, I_4, \, I_4$
fibres  $F_1, F_2, F_3, F_4$  and a section  $P_0$. 
Such an elliptic surface is unique
[MP, Th 5.4].  Write  $F_1 = C_0 + C_1, \, F_2 = D_0 + D_1, \,
F_3 = \sum_{i=0}^3 E_i$  and  $F_4 = \sum_{j=0}^3 G_i$
with  $E_i . E_{i+1} = G_j . G_{j+1} = 1$  so that
$P_0$  meets  $C_0, D_0, E_0, G_0$.
Using the height-pairing in [Sh, Th 8.6], one sees that the Mordell Weil
group of the fibration is  $\mu_2 \times \mu_4$;
after suitable relabelling, two (disjoint) sections  $P_1, P_2$
meet fibres in this way:  $P_1$  meets  $C_0, D_1, E_1, G_1$
and  $P_2$  meets  $C_1, D_0, E_1, G_3$.  Clearly
$P_1, P_2$  have order 4 and generate the Mordell Weil group
with  $P_0$  as the origin.

\par
Denote by  $g_i$  the translation automorphism given by the section
$P_i$.  One sees that  ${\widetilde X}^G$, with  
$G = \langle g_1, g_2 \rangle \cong \mu_2 \times \mu_4$,
consists of 12 points : $g_1$ (resp. $g_2$)
fixes the two nodes  $p_1, \, p_2$ (resp. $p_3, \, p_4$)
of  $F_1$ (resp. $F_2$), while  $g_2$ (resp. $g_1$) switches
$p_1, \, p_2$ (resp. $p_3, \, p_4$);
$g_2 g_1^{-1}$ (resp. $g_2  g_1$) fixes the four nodes  $p_5, \dots, p_8$
(resp. $p_9, \dots, p_{12}$)  of  $F_3$ (resp. $F_4$).
Let  $\widetilde X \rightarrow X$  
be the smooth blow-down of the eight (disjoint) sections so that  
$X = \bold P^1 \times \bold P^1$.  Then there is an induced
$G$  actions so that  $X^G$  is still a 12-point set
(the image of  ${\widetilde X}^G$).  

\par
One sees that  $Y_1 = X/H$  with  $H$ ($\cong \mu_4$) generated by 
(the image) of  $g_1$, has exactly three singular points of type  
$A_3$, $A_3$, $A_1$ (the images of  $p_1, p_2, \{p_3, p_4\}$),
$K_{Y_1}^2 = K_X^2/4 = 2$ and rank $Pic \, Y_1 = 1$.

\par
One can also verify that  $Y = X/G$  has exactly four singular points 
of type  $A_3$, $A_3$, $A_1$, $A_1$
(the images of  $\{p_1, p_2\}, \{p_3, p_4\}, \{p_5, \dots, p_8\}, \{p_9, \dots, p_{12}\}$),
rank $Pic \, Y = 1$  and  $K_Y^2 = 1$.

\par \vskip 1pc
{\bf Lemma 2.12.} {\it Let  $X$  be a del Pezzo surface with  $K_X^2 = 1$
and an effective  $\mu_5$-action such that  $\mu_5$  fixes (point wise)
a (smooth) elliptic curve  $R \in |-K_X|$.  Then modulo  $\mu_5$-equivariant 
isomorphism the pair  $(X, \mu_5)$  is equal to either  $(X_s, \mu_5)$
$(s \in \bold P^1 \setminus \{0, \pm 1\})$  or  $(X_{II}, \mu_5)$
in Example $2.4$.}

\par \vskip 1pc
{\it Proof.}  It suffices to show that the covering map
$X \rightarrow Y = X/\mu_5$  here coincides with
some one in Example 2.4.  We shall show that
the relative minimal model of the induced ellipitc fibration
on  $Y$  has a type  $II^*$  singular fibre and (only one) section.
To begin with, let  ${\widetilde X} \rightarrow X$  be the blow-up 
of the unique base point of  $|-K_X|$  with  $\widetilde E$  
the exceptional curve.
Then  ${\widetilde X}$  is a relatively minimal elliptic surface
with a section  $\widetilde E$.  Since  $X$  contains no  $(-2)$-curves for
$-K_X$  is ample, each fibre of the elliptic fibration on  ${\widetilde X}$
is irreducible.  The induced  $\mu_5$-action on  ${\widetilde X}$
fixes the proper transform of  $R$ (also denoted by  $R$,
which is a smooth fibre now) and stabilizes  $\widetilde E$.  
Clearly the rational curve  $\widetilde E$  has exactly
two  $\mu_5$-fixed points: the intersection  $\widetilde E \cap R$  and one more
point  $p$  on another fibre  $F_1$.

\par
If  $F_2$  is a singular fibre ($\ne F_1$)  then  $\{g F_2 | g \in \mu_5\}$
is a set of 5 singular fibres of the same type.  Hence
the Euler number  $12 = \chi({\widetilde X}) = \chi(F_1) + 5 t$
with  $t \ge 0$.  Thus  $\chi(F_1) = 2, 7$  because
there should be at least two singular fibres if one calculates
the Picard number in terms of contributions from fibres and
the rank of the Mordell Weil group [Sh, Cor 5.3].  

\par
Let  $\widetilde Z \rightarrow \bold P^1$  be the smooth relative minimal model of the
elliptic fibration on the quotient  $\widetilde X/\mu_5$  induced from
the one on  $\widetilde X$.  Since  $\mu_5$  acts on the base curve
of the fibration of  $\widetilde X$  as an automorphism of order 5,
we see that if  $T(\overline{F}_1)$  
is the monodromy of the fibre  $\overline{F}_1$
on  $\widetilde Z$  dominated by the fibre  $F_1$  on  $\widetilde X$, then
the monodromy  $T(F_1)$  equals  $T(\overline{F}_1)^5$.
This and  $\chi(F_1) = 2, 7$  imply that  $F_1$  is of type  $II$
and its image  $\overline{F}_1$  is of type  $II^*$ at $t = 0$
with  $t$  the inhomogeneous coordinate of the base curve 
(we arrange  $t$  this way) [BPV, Table 6 at p.159].
So  $\widetilde Z$  is a rational elliptic surface with only one section
$E$ (the image of  $\widetilde E$) and we can identify it with
either  $\widetilde Z_I$  or  $\widetilde Z_{II}$  in Example 2.4.
Then the fibre  $B$  on  $\widetilde Z$  dominated by the fibre  $R$  
is at  $t = s$  with  $s \ne 0, \pm 1$ (resp. at  $s = 1$).

\par
One sees that  $\widetilde Z \rightarrow \widetilde X/\mu_5$ 
is the the contraction of  $D_1, \, D_2$ 
in notation of Example 2.5.  Contracting further 
$\widetilde E$  on  $\widetilde X$ and (the image of)  $E$  on  
$\widetilde X/\mu_5$  to get  $X$  and  $Y = X/\mu_5$, 
we see that our  $\sigma : X \rightarrow Y$  is also a 
$\bold Z/(5)$-Galois cover totally branched at the only
singular point  $q_1$  of  $Y$  and the curve  $B$ (the image of  $R$).
It is known that such a cover is given by a relation
${\Cal O}_Y(\overline{M})^{\otimes 5} \cong {\Cal O}_Y(B)$  for
some Weil divisor  $\overline{M}$.  Since  ${\Cal O}_Y(\overline{M})
\cong {\Cal O}_Y(C)$ by {\bf 2.4.1}, our  $\sigma$  here coincides with
the one in Example 2.4.  This proves the lemma.

\par \vskip 1pc
{\bf Lemma 2.13.} {\it There is only one pair  $(X, \mu_5)$ 
of the del Pezzo surface  $X$  with   $K_X^2 = 5$  and the group
$\mu_5$  acting effectively on  $X$  modulo equivariant  $\mu_5$-isomorphism.}

\par \vskip 1pc
{\it Proof.} A degree 5 del Pezzo surface  $X$  is the blow-up of 4 points
$p_i$  on  $\bold P^2$  (no three of them are collinear), and hence there is 
only one such  $X$  modulo isomorphism (these 4 points  $p_i$  form a frame
of  $\bold P^2$, and any other frame is mapped to this by a projective
transformation).  It is known that  $Aut(X)$  is
the symmetric group  $S_5$  in 5 letters.  Since all sylow-5
groups of $S_5$  are conjugate to each other, the lemma follows.

\par \vskip 1pc
{\bf Lemma 2.14.} {\it Modulo equivariant  $\mu_6$-isomorphism,
there is only one pair  $(X, \mu_6)$  of the $($unique$)$ del Pezzo surface  $X$  
with   $K_X^2 = 6$  and the group  $\mu_6$  acting effectively on  $X$
such that  $X^{\mu_6}$  is a finite set and  $Y = X/\mu_6$  has exactly
$3$  singularities of type  $\frac{1}{2}(1,1), \, \frac{1}{3}(1,2), \,
\frac{1}{6}(1,5)$  as all of its singularities.}

\par \vskip 1pc
{\it Proof.} Note that  $Y$  is a Gorenstein log del Pezzo surface 
with  $3$  singularities of type  $A_1, A_2, A_5$  in other notation.
In view of Lemma 4.4, it is enough to show that there is only one such
$Y$  modulo isomorphism.
Let  $f : Z \rightarrow Y$  be the minimal resolution.
Then  $Z$  is an almost del Pezzo surface with  $K_Z^2 = 1$
so that  $|-K_Z|$  has exactly one base point, $\dim |-K_Z| = K_Z^2 = 1$
(Riemann-Roch and Kawamata-Viehweg vanishing) and a general member
of  $|-K_Z|$  is irreducible [Dem, Th 1, page 39].  Let  $\widetilde Z \rightarrow Z$
be the blow-up of the unique base point of  $|-K_Z|$  with  $P_0$
the exceptional curve.  Then  $\widetilde Z$  is a relatively
minimal elliptic surface with  $P_0$  as a section.  One sees that
the inverse of  $\text{\rm Sing} \, Y$  is contained in three different
fibres  $F_1, F_2, F_3$  of types  $I_2, I_3, I_6$ [MP, Theorem 4.1].
Now the uniqueness of  $Y$  follows from 
the uniqueness of such elliptic surface [MP, Th 5.4] and
also the uniqueness of the pair  $(\widetilde Z, P_0)$
modulo translation automorphism.  This proves the lemma.

\par \vskip 1pc
{\bf Lemma 2.15.} (1) {\it Modulo equivariant  $H$-isomorphism ($H \cong \mu_4$)
there is only one pair  $(X, H)$  with  $X = \bold P^1 \times \bold P^1$
and the group  $H$  acting effectively on  $X$
such that  $X^H$  is a finite set and  $Y_1 = X/H$  has exactly
$3$  singularities of type  $\frac{1}{2}(1,1), \, \frac{1}{4}(1,3), \,
\frac{1}{4}(1,3)$  as all of its singularities and rank $Pic \, Y_1 = 1$.}

\par
(2) {\it Modulo equivariant  $G$-isomorphism ($G \cong \mu_4 \times \mu_2$),
there is only one pair  $(X, G)$  with  $X = \bold P^1 \times \bold P^1$
and the group  $G$  acting effectively on  $X$
such that  $X^G$  is a finite set and  $Y = X/G$  has exactly
$4$  singularities of type  $\frac{1}{2}(1,1), \, \frac{1}{2}(1,1), \, 
\frac{1}{4}(1,3), \, \frac{1}{4}(1,3)$  as all of its singularities.}

\par
(3) {\it There is a subgroup  $H_1$  of  $G$  such that
$(X, H_1) = (X, H)$  modulo  $\mu_4$-equivariant isomorphism
(identify  $H = H_1 = \mu_4$) and hence  $Y = X/G = (X/H_1)/\overline{G}
= Y_1/\overline{G}$, where  $\overline{G} = G/H_1 \cong \mu_2$.}

\par \vskip 1pc
{\it Proof.} (1) Note that  $Y_1$  is a Gorenstein log del Pezzo surface
of Picard number 1 and with singularities of type  $A_1, A_3, A_3$.
It suffices to show such  $Y_1$  is unique (Lemma 4.4).
Let  $f : Z_1 \rightarrow Y_1$  be the minimal resolution.
Then  $f^{-1}(\text{\rm Sing} \, Y_1)$  is a disjoint union of 
linear chains of  $(-2)$-curves of length $1, 3, 3$.
As in Lemma 2.14, $Z_1$  is an almost del Pezzo surface with  $K_{Z_1}^2 = 2$,
$\dim |-K_{Z_1}| = 2$  and a general member  $A$  of  $|-K_{Z_1}|$  smooth
irreducible [Dem, Th 1, page 39].  Pick up any  $(-1)$-curve  $D_4'$  on  $Z_1$
such that  $D_4'+f^{-1}(\text{\rm Sing} \, Y_1)$  is a disjoint union
of two linear chains of length $1, 7$
($D_4'$  connects the two length-3 chains).  One can find such
$D_4'$  by playing with  $\bold P^1$-fibrations, or from
[Zh1, Lemmas 3.5, 4.2, 4.3] we see that  $(Z_1, f^{-1}(\text{\rm Sing} \, Y_1))$
fits Case (9) in Lemma 4.2 there and the picture at [Zh1, p. 454], 
and we just let  $D_4' = E_2$  in notation there.  

\par
Let  $Z_2 \rightarrow Z_1$  be the blow-up of the point
$A \cap D_4'$  (with  $A$  fixed for the time being)  with  $P_2'$
the exceptional curve.
Then  $Z_2$  is again an almost del Pezzo surface with  
$K_{Z_2}^2 = 1$.  As in Lemma 2.14, let  
$Z_3 \rightarrow Z_2$  be the blow-up of
the only base point of  $|-K_{Z_2}|$  with  $P_0$  the exceptional curve.
Then  $Z_3$  is a relatively minimal elliptic surface
so that the proper inverse of  $D_4'+f^{-1}(\text{\rm Sing} \, Y_1)$  
are contained in two different fibres  $F_1, \, F_2$  of types  $I_2, \, I_8$
[MP, Th 4.1].
Now the uniqueness of  $Y$  follows from 
the uniqueness of such elliptic surface [MP, Th 5.4] and
also the uniqueness of the triplet  $(Z_3; P_0, P_2)$
modulo translation automorphism, noting that the Mordell
Weil group of the fibration is  $\bold Z/(4)$  and when
we choose  $P_0$  as the origin then  $P_2$  is the unique
element of order 2.  This proves (1).

\par
(2) Note that  $Y$  is a Gorenstein log del Pezzo surface 
with  $4$  singularities of type  $A_1, A_1, A_3, A_3$.
As in Lemma 2.14,
let  $f : Z \rightarrow Y$  be the minimal resolution
and let  $\widetilde Z \rightarrow Z$
be the blow-up of the unique base point of  $|-K_Z|$  with  $P_0$
the exceptional curve. 
Now the uniqueness of  $Y$  follows from 
the uniqueness of such elliptic surface [MP, Th 5.4] and
also the uniqueness of the pair  $(\widetilde Z, P_0)$
modulo translation automorphism.  This proves (2).

\par
(3) is shown in Example 2.11.

\par \vskip 2pc
{\bf \S 3. Case : the invariant sublattice is of rank 1}

\par \vskip 1pc
In this section, we consider the case ($\rho(Y) = $) 
rank $(Pic X)^{\mu_p} = 1$.  We first treat the case  $r(Y) < 1$.

\par \vskip 1pc
{\bf Proposition 3.1.} {\it Suppose that the quotient surface satisfies 
$\rho(Y) = 1$  and  $r = r(Y) < 1$.  Then  $p \ge 5$, $(X, \mu_p)$
equals a pair in Example $2.1b$, and  $r = 1/p$  or  $r = 3/p$.}

\par \vskip 1pc
{\it Proof.} By Lemmas 1.2 and 1.8, $X^{\mu_p}$  is a finite set
$\{p_1, \cdots, p_s\}$  with  $s \ge 1$.  Set  $q_i = \sigma(p_i)$  so that
Sing $Y = \{q_1, \dots, q_s\}$.
Note that if  $X = \bold P^2$  and  $s = 3$  then  $p \ge 5$
because  $r < 1$ (see Example 2.1).  Write  $r = m/p$
with an integer  $1 \le m \le p-1$.

\par
We shall frequently apply Lemma 1.6 to the extent
that  $9 - (p-1)(3-s) = K_X^2$, which is between 1 and 9 because  $X$  is del Pezzo
(Lemma 1.7).  Since  $2g(P) - 2 = P . (P + K_Y) = (p-m)P^2/p$  is an integer (Lemma 1.7), 
$p | P^2$.  One has also  $K_X^2 = (\sigma^*K_Y)^2 = p K_Y^2 = m^2 P^2/p \ge m^2$.

\par
If  $m \ge 3$, then  $K_X^2 = 9$, $r = 3/p$, 
$P^2 = p$  and  $s = 3$;
so the Proposition is true.
If  $m = 2$  and  $P^2/p \ge 2$, then  $K_X^2 = 8$, $r = 2/p$
and  $(p, s) = (2, 2)$, which leads to that  $r = 1$, a contradiction.
If  $m = 2$  and  $P^2 = p$, then
$K_X^2 = 4$, which leads to  $4 = 9 - (p-1)(3-s)$, a contradiction.

\par
We now assume that  $r = 1/p$.  
If  $s = 3$, then  $K_X^2 = 9$, $P^2 = 9p$; 
so the Proposition is true.

\par
Suppose that  $s = 1$.  Then  $K_X^2 = 9 - 2(p-1)$, whence  $p = 2, 3, 5$.
In notation of Lemma 1.6, one has  $k_u = s = 1$  for some  $1 \le u \le p-1$
and  $1 = k_u (1/(p-1))\sum_{\zeta} 1/(1 - \zeta)(1 - \zeta^u)$,
where  $\zeta$  runs over the set of all primitive  $p$-th root of 1.
This is impossible by the calculation of the right hand side in Lemma 1.6.

\par
Suppose that  $s = 2$.  Then  $K_X^2 = 10 - p$  and  $p = 2, 3, 5, 7$.
Since  $\mu_7$  can not act on a cubic del Pezzo
(see [Man2, Table 1, p.176]), $p \ne 7$.
Applying Lemma 1.6, we can show that
$p = 5, k_4 = s = 2$.  But then the quotient surface  $Y$  is Gorenstein
and hence  $r(Y) \ge 1$, a contradiction.
This completes the proof of the Proposition.

\par \vskip 1pc
Next we consider the case  $r(Y) = 1$.

\par \vskip 1pc
{\bf Proposition 3.2.} {\it Suppose that the quotient surface satisfies
$\rho(Y) = 1$  and  $r(Y) = 1$.  Then  $(X, \mu_p)$  is equal to one of the pairs 
in Examples $2.1b$ (with  $p = 3$), $2.4$  and  $2.5$.}

\par \vskip 1pc
{\it Proof.} We note that the minimal resolution  $Z$  of  $Y$  is neither  $\bold P^2$
nor a Hirzebruch surface  $\bold F_e$, for otherwise, either  $Y = Z = \bold P^2$,
or  $Y = \overline{\bold F}_e$ ($e \ge 2$)  because  $\rho(Y) = 1$, which 
would lead to  $r(Y) > 1$, a contradiction (Remark 1.5);
in particular, $K_Z^2 \le 7$.  If  $K_X^2 \ge 8$, then
$X = \bold P^2$  or  $X = \bold F_e$ ($e = 0, 1$) because  $X$  is del Pezzo (Lemma 1.7).
If  $X = \bold F_1$  then  $\mu_p$  stabilizes the  $(-1)$-section
and the divisor class of a fibre, which contradicts that  $\rho(Y) = 1$ (Lemma 1.2).
By the same reasonning, if  $X = \bold F_0 = \bold P^1 \times \bold P^1$, 
then  $p = 2$  and  $\mu_p$  switchs the two fibrations, but then
$Y = \bold P^2$  and  $r = 3 > 1$, a contradiction (see Example 2.6).
If  $X = \bold P^2$, then the Proposition is true by Lemma 1.9 (see Example 2.1).
So we may assume the following, noting that  $f^*K_Y = K_Z$ (Lemmas 1.3 and 1.9) :

\par \vskip 1pc
{\bf Condition 3.2.1.} $1 \le K_X^2 \le 7$  and  $1 \le K_Y^2 = K_Z^2 \le 7$.

\par \vskip 1pc
Consider first the case  $X^{\mu_p}$  contains an irreducible curve  $R$.
We shall show that this fits Example 2.4.
By Lemma 1.7, $X^{\mu}$  is a union of the irreducible (smooth) 
curve  $R$  and  $s$  points  $p_i$.  As in Lemma 1.8, writing  $\sigma(R) = B$  
and  $B = b P = b(-K_Y)$  with  $b$  a positive integer (Lemma 1.7), we get
$-K_X = (m/n) \sigma^* P$, where  $m/n = [p - b(p-1)]/p$  with coprime
positive integers  $m, n$.  So  $1 \le b = [1 - (m/n)][p/(p-1)] < p/(p-1) \le 2$.
Thus  $b = 1$  and  $m/n = 1/p$.
So  $B$ ($\sim -K_Y$) and its (reduced) preimage  $R$  are 
isomorphic elliptic curves.
Now  $-K_X = (1/p) \sigma^* P \sim R$  and  $K_X^2 = P^2/p = K_Y^2/p$.
By Lemma 1.6, one has also  $K_X^2 = 9 - (p-1)(3-s)$.

\par
Since each  $q_i \in \text{\rm Sing} \,  Y$ ($1 \le i \le s$) is Du Val, 
$D_i := f^{-1}(q_i)$  is a 
chain of  $p-1$  of  $(-2)$-curves.  So  $\rho(Z) = \rho(Y) + s(p-1) = 1 + s(p-1)$.
Thus  $K_Y^2 = K_Z^2 = 10 - \rho(Z) = 9 - s(p-1)$, which is an integer
between  $1$  and  $7$ (cf. 3.2.1).
So  $2 \le s(p-1) \le 8$.
Solving  $p[9 - (p-1)(3-s)] = pK_X^2 = K_Y^2 = 9 - s(p-1)$,
one obtains  $s = 3 - 12/(p+1)$.  So only  $(p, s) = (5, 1)$
is possible.  Our pair here is equal to the pair in Example 2.4
modulo  $\mu_5$-equivariant isomorphism (Lemma 2.12).

\par
Next we consider the case where  $X^{\mu_p}$  is a finite set
$\{p_1, \dots, p_s\}$  with  $s \ge 1$ (Lemma 1.2).
Then  $-K_X = \sigma^*K_Y$  and  $K_X^2 = p K_Y^2$.
As above, one obtains  $9 - (p-1)(3-s) = K_X^2 = p K_Y^2 =
p[9 - s(p-1)]$  with  $2 \le (p-1)(3-s) \le 8$ (cf. 3.2.1),
$s = 12/(p+1)$, and  $(p, s) = (5, 2)$.  Our pair here is now
equal to the pair in Example 2.5 modulo  $\mu_5$-equivariant isomorphism 
(Lemma 2.13).  This completes the proof of the Proposition.

\par \vskip 1pc
Now we treat the case where  $r(Y) > 1$. 

\par \vskip 1pc
{\bf Proposition 3.3.} {\it Suppose that the quotient surface satisfies
$\rho(Y) = 1$  and  $r(Y) > 1$.  Then  $(X, \mu_p)$
is equal to one of the pairs in Examples $2.1a, 2.2, 2.3, 2.6, 2.7, 2.8$.}

\par \vskip 1pc
{\it Proof.} By Lemma 1.10, (a general member) $P$  is a smooth rational curve
away from the singular locus of  $Y$.  Let  $P' = f^*P$  and  $m = P^2$.
Applying the cohomology exact squence arising from the exact sequence:
$$0 \rightarrow {\Cal O}_Z \rightarrow {\Cal O}(P') 
\rightarrow {\Cal O}_{P'}(m) \rightarrow 0,$$
we obtain  $h^0(Z, P') = (P')^2+2$.
As long as  $L$  is a smooth rational curve with  $L^2 \ge 0$  
on a smooth rational surface,
one always has  $h^0(L) = L^2 + 2$; thus by induction on  $L^2$
(to reduce to  $L^2 = 0$ case) one can deduce that  $Bs|L| = \emptyset$.

\par
So the linear system  $|f^*P|$  gives rise to a well-defined morphism
$\Phi : Z \rightarrow \bold P^{m+1}$, with the image  $W$  a non-degenerate
surface (noting that  $P^2 > 0$)  and hence  deg $W \ge m$.
On the other hand, $m = (P')^2 = (\text{\rm deg} \Phi) (\text{\rm deg} W)$.  Thus
$\Phi$  is a birational morphism onto a degree  $m$  surface
in  $\bold P^{m+1}$.  Clearly, $\Phi$  factors as  $\varphi \circ f$,
with a birational morphism  $\varphi : Y \rightarrow W$
which is given by the linear system  $|P|$.
Since  $P$  is ample, $\varphi$  is an isomorphism ($W$  is normal; see below).
So we can identify  $Y = W$.

\par
Non-degenerate surfaces  $W$  of degree  $m$  in  $\bold P^{m+1}$
are well classified (cf. [Nag]).
$W$  is either  $\bold P^2$ ($m = 1$), or the Veronese embedding of
$\bold P^2$  in  $\bold P^5$ ($m = 4$), or the embedding (with the negative 
section  $C$  contracted if  $a = n = m$)
of the Hirzebruch surface  $\bold F_n$  in  $\bold P^{m+1}$ 
by the linear system  $|C + aF|$, where  $m = 2a -n, a \ge n$, 
$C$  is the section with  $C^2 = -n$  and  $F$  a fibre.
If  $W$  is smooth then  $Z = Y = W \cong \bold P^2$  because
$\rho(Y) = 1$.  If  $W$  is singular, then  $a = n = m \ge 2$, $Z = \bold F_m$
and  $Y = \overline{\bold F}_m$, the projective cone (see Remark 1.5).
Clearly, $m = p$  and the only singularity
in  $\overline{\bold F}_p$  is of type  $\frac{1}{p}(1,1)$.

\par
Suppose that  $Y = \bold P^2$.  Then  $X^{\mu_p}$  is a single
smooth curve  $R$ (Lemma 1.7).  Let  $d$  be the degree in  $Y$  of  $B = \sigma(R)$.
Then the  $\bold Z/(p)$-Galois cover  $\sigma : X \rightarrow Y$
is given by a relation  $B \sim p F$.  Set 
$d = \text{\rm deg}(B) = p \text{\rm deg}(F)$.
Now  $K_X = \sigma^*(K_Y + (p-1)F) = [(p-1)(deg F) - 3] \sigma^* P$,
where  $P$  is a line.  Since  $-K_X$  is ample, $(p, \text{\rm deg} F) =
(2, 1), (2, 2), (3, 1)$.  Thus  $(X, \mu_p)$  is as in
Examples 2.6, 2.7 and 2.2.

\par
Suppose that  $Y = {\overline {\bold F}}_p$.  If  $X^{\mu_p}$
contains no curve, then it is a single point  $p_1$  and
$\sigma$  is unramified over  $Y - \{q_1\}$, where
$q_1 = \sigma(p_1)$  is the vertex of the cone  $Y$;
this is impossible because  $Y - \{q_1\}$  is simply connected.
Write  $B = bP$.  This  $P$  is the generator
of  $Pic \, Y$, is the hyperplane of  $Y \subseteq \bold P^{p+1}$
and satisfies  $P^2 = p$ ($P = \overline H$  in notation of Remark 1.5). 
In the present case, $b \ge 1$  is an integer
for  $B \in Pic \, Y$.  As in Lemma 1.8, $-K_X = [(rp - b(p-1)]/p \sigma^*P
= [(p+2) - b(p-1)]/p \sigma^*P$, noting that
$r = (p+2)/p$ (Remark 1.5).  Since  $-K_X$  is ample,
either  $b = 1$, or  $(b, p) = (2, 2), (2, 3), (3, 2)$.

\par
Note that the  $\bold Z/(p)$-Galois cover  $\sigma : X \rightarrow Y$  
is totally branched along the smooth curve  $B$  and
the vertex  $q_1$ ($\notin B$).
If  $b = 1$, then one sees easily that  $(X, \mu_p)$  is equal
to a pair in Example 2.1a.

\par
If  $(b, p) = (2, 2)$, one can verify that the  $\sigma^{-1}(q_1)$
has to split into two singularities of the same type, i.e.,
$\sigma$  is not branched at  $q_1$, a contradiction.
If  $(b, p) = (2, 3)$  or  $(3, 2)$  then  $(X, \mu_p)$  
is as in Example 2.3 or 2.8.  This proves Proposition 3.3.

\par \vskip 1pc
{\bf \S 4.  The proofs of Theorems and Corollaries}

\par \vskip 1pc
Let  $(X, G)$  be a pair with  $X$  a smooth projective rational surface
and  $G$  a finite group acting effectively on  $X$.
We follow the approach in [BB] using the Mori theory.
The cone theory [Mor, Theorems 1.5 and 2.1] implies
the decomposition of the closed cone of effective cycles
with coefficients in  $\bold R$  and modulo numerical equivalence:
$$\overline{NE}(X) = \overline{NE}(X)_{K_X \ge 0}+ \sum_{C \in \Cal E} \bold R_+ [C],$$
where  $\Cal E$  is a countable set of smooth rational curves  $C$
satisfying  $C^2 = -1, 0, 1$.  Passing to the  $G$-invariant part,
we get [Mor, Proposition 2.6]:
$$\overline{NE}(X)^{G} = \overline{NE}(X)^{G}_{K_X \ge 0} 
+ \sum_{C \in \Cal F} \bold R_+ [\sum_{g \in G} g C],$$
where  $\Cal F$  is the subset of curves  $C$  in  $\Cal E$  such that
$\bold R_+ [\sum_{g \in G} g C]$  is an extremal ray in the
$G$-invariant cone  $\overline{NE}(X)^{G}$.

\par
For a curve  $C$  on  $X$, denote by  $G_c$  the maximum subgroup of  
$G$  stabilizing  $C$  and let  $k_c$  be the index  $|G : G_c|$  and  
$\{g_i G_c | 1 \le i \le k_c\}$  (with  $g_1 =$ id) the  $k_c$  cosets.
For the lemma below, we are essentially proceeding along the idea in 
[Mor, Theorem (2.7)].

\par \vskip 1pc
{\bf Lemma 4.1.} {\it Assume that  $(X, G)$  is a minimal pair
such that $\text{\rm rank}(Pic \, X)^G \ge 2$.
Then there is a  $G$-stable conic fibration $\varphi : X \rightarrow \bold P^1$
with a smooth rational curve as its general fibre, such that
every singular fibre is a linear chains of two  $(-1)$-curves.
If  $\varphi$  is not smooth, i.e. $X$  is not a Hirzebruch surface
then  $|G|$  is even.}

\par \vskip 1pc
{\it Proof.} Since  $X$  is rational,
$K_X$  has negative intersection with a curve  $E$  and hence
with the  $G$-stable effective 1-cycle  $\sum_{g \in G} gE$.
So the cone  $\overline{NE}(X)^{G}$  has an extremal ray
$\bold R_+ [L]$  where  $L = \sum_{g \in G} g C$
with a smooth rational curve  $C$.

\par
Note that  $L_{\text{\rm red}} = \sum_{i=1}^k g_i C$
and  $L = |G_c| L_{\text{\rm red}}$, where $k = k_c$.
First  $L^2 \le 0$, for otherwise  $L$  would belong to the
interior of  $\overline{NE}(X)^G$ [Mor, Lemma (2.5)]
and  $L$  could not be extremal; here we use the fact that
$\text{\rm rank}(Pic \, X)^G \ge 2$.  Since  $L^2 \le 0$, we have  $C^2 \le 0$;
and if  $C^2 = 0$  then  $L_{\text{\rm red}}$  
is a disjoint union of  $k$  smooth rational curves of self intersection 0
and  $\varphi := \Phi_{|C|}$  is a  $G$-stable  $\bold P^1$-fibration.

\par
Consider the case  $C^2 \le -1$.
Then  $C$  is a  $(-1)$-curve  for  $C . K_X < 0$.
If  $k = 1$, then  $C$  is  $G$-stable, a contradiction to
the minimality of the pair.  So  $k \ge 2$.
Now  $0 \ge L^2 = |G| (C . L) = |G||G_c| (C . L_{\text{\rm red}})
= |G||G_c| (-1 + \sum_{i \ge 2} C . g_i C)$.
If  $\sum_{i \ge 2} C . g_i C = 0$  then
$L_{\text{\rm red}}$  is a disjoint union of  $k$  of  $(-1)$-curves
which contradicts the minimality of the pair.  So
we may arrange as  $C . g_2 C = 1$  and  $C . g_i C = 0$ ($i \ge 3$).
Since  $g_2^{-1} C . C = 1$, one has  $g_2^{-1} C = g_2 C$
and  $g_2^2 \in G_c$  because  $C$  meets only
one component  $g_2 C$  among  $g_i C$ ($2 \le i \le k$).

\par
If  $k \ge 3$, we see that  $g_3 (C + g_2 C)$  is a linear
chain of two intersecting  $(-1)$-curves disjoint from
$C + g_2 C$.  We can easily arrange  $L_{\text{\rm red}}$
as a disjoint unions of pairs of intersecting  $(-1)$-curves
$g_i (C + g_2 C)$ ($1 \le i \le k/2$)
so that an arbitrary element of  $G$  either stabilizes each of
two components in  $C + g_2 C$, or switch them or maps them to
some  $g_i (C + g_2 C)$.  Thus  $\varphi = \Phi_{|C + g_2 C|}$
is a  $G$-stable  $\bold P^1$-fibration.  Note that  $2$  divides
$k = k_c = |G : G_c|$  and hence  $|G|$.

\par
To finish the proof, we still need to determine the singular fibres
of the  $G$-stable conic fibration  $\varphi$.
Take any  $(-1)$-curve  $E$  in a fibre of  $\varphi$
and set  $L = \sum_{g \in G} g E$.  Then  $L^2 \le 0$
because  $L$  is supported by fibres and hence negative semi-definite [Re2, A.7].
Now the same argument above will imply the lemma.

\par \vskip 1pc
{\bf Corollary 4.2.} {\it Assume that  $p$  is an odd prime number.
Let  $(X, \mu_p)$  be a minimal pair with
$X$  a smooth projective rational surface and
$\mu_p$  acting effectivly on  $X$.
Suppose that  $(Pic X)^{\mu_p}$  has rank  $\ge 2$.
Then  $X$  is a Hirzebruch surface  $\bold F_e$ ($e \ne 1$)
and every ruling (there are two only when  $e = 0$) is  $\mu_p$-stable.}

\par \vskip 1pc
{\it Proof.} By Lemma 4.1, $X = \bold F_e$.
If  $e = 1$, then the unique  $(-1)$-curve
would be  $\mu_p$-stable and we reach a contradiction to the
minimality assumption.  The rest is clear.

\par \vskip 1pc
Before we proceed to prove theorems, we need two results.

\par \vskip 1pc
{\bf Lemma 4.3.} {\it Suppose that the group  $\mu_p$  
of prime order  $p$  acts effectively on the Hirzebruch surface  $X = \bold F_e$
and stabilizes its fixed ruling  $\varphi$.  Then  $X^{\mu_p}$  is one of the
following, where  $\mu_p$  stabilizes exactly two out of all fibres in 
the first three cases,}

\par 
(1) {\it a union of two fibres,}

\par
(2) {\it a union of a fibre and two points in another fibre},

\par
(3) {\it four points in two distinct fibres, and}

\par
(4) {\it a disjoint union of two sections (one of which is the
$(-e)$-section if  $e > 0$).}

\par \vskip 1pc
{\it Proof.} This result must be well known but we do not
have reference.  Suppose that  $X^{\mu_p}$  is contained in fibres.
Note that there is an induced  $\mu_p$-action on the base
rational curve of the ruling.  So either  $\mu_p$  stabilizes
exactly two fibres (then Case (1), (2) or (3) of the lemma
occurs), or  $\mu_p$  stabilizes all fibres. 
If the second situation happens, then  $X = \bold F_0 = \bold P^1 \times \bold P^1$,
for otherwise the  $(-e)$-section is  $\mu_p$-stable only and has
exactly two  $\mu_p$-fixed points so that
$\mu_p$  would stabilize only the two fibres containing these two points;
on the other hand, $\mu_p$  stabilizes also the second ruling
as well as its fibre  $F_2$  through a point in  $X^{\mu_p}$
(which is non-empty; see Lemma 1.2), so that  $F_2$  is
a section of the first ruling and has exactly two  $\mu_p$-fixed
points, a contradiction again. 

\par
Next we consider the case where  $X^{\mu}$  contains a (multi-)section.
Then each fibre is  $\mu_p$-stable.
Thus either  $X^{\mu}$  is the union of two disjoint sections
(one of which is the  $(-e)$-curve if  $e > 0$) 
so that Case (4) of the lemma occurs,
or  $X^{\mu}$  is the union of a double section  $D$  and a few points
(a general fibre of the ruling has exactly two  $\mu_p$-fixed points).
In the second case, $X = \bold F_0$  
and  $\mu_p$  stabilizes also the second ruling.
Since  $D$  intersects all fibres of both rulings, an arbitrary
fibre of any ruling is  $\mu_p$-stable, whence the diagonal
of  $X$ (a section of both rulings) is also contained in  $X^{\mu_p}$, 
a contradiction to the assumption that  $X^{\mu_p}$  is the
union of a double section  $D$  of the first ruling  $\varphi$
and a few points.  This proves the lemma.

\par \vskip 1pc
{\bf Lemma 4.4.} {\it For  $i = 1, 2$, let  $(X_i, G)$  be a 
pair of a simply connected smooth algebraic surface (e.g. a rational surface)
and a finite group  $G$  acting effectively on  $X_i$  such that
$X_i^{G}$  is a finite set.  Let  $Y_i = X_i/G$  and  
$Y_i^0 = Y_i - \text{\rm Sing} \, Y_i$.  Then we have:}

\par
(1) {\it One has  $\pi_1(Y_i^0) = G_i$  and the quotient map  
$\sigma_i : X_i \rightarrow Y_i$  is the completion of the universal
covering map  $U_i^0 \rightarrow Y_i^0$; in other words,
$X_i$  is the normalization of  $Y_i$  in the function field
${\bold C}(U_i^0)$.}

\par
(2) {\it Two pairs  $(X_i, G)$  are equal modulo  $G$-equivariant isomorphism
if and only if the  $Y_i$  are isomorphic to each other.}

\par \vskip 1pc
{\it Proof.} Since  $X_i$  with a few points removed, is 
still simply connected, (1) follows.
(2) is a consequence of (1).

\par \vskip 1pc
Now we prove Theorem 1.  Theorem 1 (II) is a consequence of
Propositions 3.1, 3.2 and 3.3.  For Theorem 1 (I), in view of Corollary 4.2,
we only need to show that every pair  $(\bold F_e, \mu_p)$  
is  $\mu_p$-birationally equivariant to a pair  $(\bold P^2, \mu_p)$
in Example 2.1.  This can be proved as in [BB, (2.5)].  For readers' convenience,
we give a sketch here.  Let  $\varphi$  be as in Corollary 4.2.

\par \vskip 1pc
{\bf Case 4.5}: $\mu_p$  stabilizes each fibre.
When  $e > 0$  the unique  $(-e)$-section is  $\mu_p$-fixed
and  $X^{\mu}$  contains one more disjoint section.
We blow up a point  $p_1$  on the second (positive) section
and blow down the proper transform of the fibre containing  $p_1$.
Then we get a  $\mu_p$-birational equivariance between
our original pair  $(\bold F_e, \mu_p)$  and a new pair
$(\bold F_{e-1}, \mu_p)$.  Inductively we reduce to
the case  $e = 1$  and further blow down the  $\mu$-stable  $(-1)$-curve
on  $\bold F_1$  to proceed  $\mu_p$-birationally equivariantly to
a pair  $(\bold P^2, \mu_p)$  in Example 2.1.

\par \vskip 1pc
{\bf Case 4.6}: $\mu_p$  acts non-trivially on the set of fibres
(and hence on the base rational curve of ruling).
Then there are exactly two  $\mu_p$-stable fibres
(lying over two  $\mu_p$-fixed points of the base rational curve;
see Lemma 4.3).
Each stable fibre contains at least two  $\mu_p$-fixed points.
We blow up the one not lying on the  $(-e)$-curve and then
blow down the proper transform of the fibre; we reduce
to a pair  $(\bold F_{e-1}, \mu_p)$.  The rest is the same as
in Case 4.5.

\par
For both Cases (4.5) and (4.6), when  $e = 0$, the argument is similar (see [BB, (2.5)]).
This completes the proof of Theorem 1.

\par \vskip 1pc
For Corollary 2, we first proceed  $\mu_p$-birationally equivariantly to
a minimal pair and then apply Theorem 1.  

\par \vskip 1pc
For Corollary 3, it suffices to consider minimal pairs.
Indeed, start with a pair  $(X, \mu_p)$  and
let  $(X_{\min}, \mu_p)$  be a minimal pair
with a  $\mu_p$-equivariant birational {\it morphism}
$\tau : X \rightarrow X_{\min}$; then 
$\tau$  induces a birational morphism  $\tau_y : Y = X/\mu_p \rightarrow
Y_{\min} = X_{\min}/\mu_p$; the images of $ \text{\rm Sing} \, Y$
and the  $\tau_y$-exceptional divisor form a finite subset  $\Sigma$  of  $Y_{\min}$,
and  $Y_{\min}^0 \setminus \Sigma$  can be regarded as a Zariski-open
subset of  $Y^0$, whence we have a surjective homomorphism  
$\pi_1(Y_{\min}^0) = \pi_1(Y_{\min}^0 \setminus \Sigma) \rightarrow
\pi_1(Y^0)$.  

\par \vskip 1pc
{\bf Remark 4.7.} In particular, if  $X^{\mu_p}$
contains a curve  $R$  with  $R^2 \ge 0$  or  $g(R) \ge 1$
then the image on  $X_{\min}$  of  $R$  is still a {\it curve}
in  $X_{\min}^{\mu_p}$, whence  $\pi_1(Y^0) = \pi_1(Y_{\min}^0) = (1)$
by the statement for minimal pairs.

\par \vskip 1pc
We now prove Corollary 3 for a minimal pair  $(X, \mu_p)$.
If the lattice  $(Pic X)^{\mu_p}$  has rank 1, then
Corollary 3 is true by Table 1 in Theorem 1.
Suppose that this lattice has rank  $\ge 2$.  Then
$X$  is a Hirzebruch surface  $\bold F_e$
and the fixed ruling  $\varphi : X \rightarrow \bold P^1$ 
is  $\mu_p$-stable (Corollary 4.2).  If  $X^{\mu_p}$
is a finite set then  $X \setminus X^{\mu_p}$  is simply connected
and equals the universal cover of  $Y^0$, whence Corollary 3 is true.
If  $X^{\mu_p}$  is a (disjoint) union of smooth curves
then  $Y = X/\mu_p$  is smooth rational and hence
$Y^0 = Y$  is simply connected.  

\par
It remains to consider Lemma 4.3,Case(2).
Then  $Y$  is rational with two singular points
$q_i$ (images of two isolated  $\mu_p$-fixed points) 
and a ruling  $Y \rightarrow \bold P^1$
(induced from the one on  $X$) such that both  $q_i$  
are on the same fibre  $F_1$.
Thus  $Y \setminus F_1$  is a  $\bold P^1$-bundle over
the affine line  $\bold A^1$  and hence simply connected.
Now the inclusion  $Y \setminus F_1 \subseteq Y^0$
induces a surjective map $(1)=\pi_1(Y\setminus F_1)
\rightarrow\pi_1(Y^0)$, whence $\pi_1(Y^0) = (1)$.
This proves Corollary 3. 

\par \vskip 1pc
We prove Theorem 4. (I) follows from Lemma 4.1.
Next we do (II) and so assume $(Pic \, X)^G$  has rank 1.
Then by Lemma 1.7, $X$  is del Pezzo and  $Y$  is singular del Pezzo.
Now the main theorem in [GZ1,2] shows that  $\pi_1(Y^0)$
is finite (see [FKL] for a differential geometric proof
and also [MS] for a new proof).  Since the rational surface  $X$  with
a few points removed, is still simply connected, it is clear that
$\pi_1(Y^0) = G$  and  $\sigma : X \rightarrow Y$  is the completion
of the universal covering map of  $Y^0$  provided that
$X^G$  is a finite set (Lemma 4.4).
Theorem 4 (II) (2) follows from the first half of the arguments
in Proposition 3.3.

\par
We now prove Theorem 4 (II) (1).  So assume  $r(Y) = 1$  and  $X^G$  is finite.
Then  $Y$  is a Gorenstein log del Pezzo surface of Picard number 1
so that  $G = \pi_1(Y^0)$  and  $\sigma : X \rightarrow Y$  is the completion 
of the universal covering map of  $Y^0$ (see Lemma 4.4).  Such  $Y$  is classified
in [Fur, Th 2] or [MZ, Lemma 6]; see also [BBD, page 593], [Ura].  
Since our  $X$  is smooth, by [MZ, Table 1, p. 71], $(X, G)$  fits
one of the rows in Table 2, but the column on  $X^G; G = \langle g_1, \dots \rangle$
is still to be verified.  For rows 2, 3, 4, 6, this is done in the
examples in \S 2 since we have the uniqueness by Lemmas 2.13-15.
For rows 1 and 5, the generator(s) of  $G$  can be easily diagonalized 
as in Table 2 (see Examples 2.1 and 2.10).
This completes the proof of Theorem 4.

\par \vskip 1pc
{\bf Final Remark 4.8.} In [Zh2, Appendix], 
there are examples of non-abelian finite group
acting effectively on  $X = \bold P^2$, such that the fixed locus  $X^G =
\{x \in X | gx = x \,\, \text{\rm for some} \,\, 1 \ne g \in G\}$
is a finite set.  For instance,
the non-abelian group of order  $21$  can act on  $\bold P^2$
this way.  Also shown are examples with  $\bold P^2$  replaced by
smooth del Pezzo surfaces or projective cones  $\overline{\bold F}_e$;
e.g., $X = \bold P^1 \times \bold P^1$  with an effective action by 
a non-abelian group  $G$  of order  $16$  or  $20$.

\par \vskip 2pc \hskip 4pc
{\bf Appendix. Weyl group based proof of Table 1}

\par \vskip 1.5pc \hskip 12pc 
I. Dolgachev

\par \vskip 2pc
Let  $(X, \mu_p)$  be a minimal pair of a smooth rational projective
surface and the group  $\mu_p$  of odd prime order  $p$
acting effectively on  $X$.  Write  $\mu_p = \langle g \rangle$.
Assume that  $X \ne \bold P^2$  and  $(Pic \, X)^{\mu_p}$
has rank 1.  In this section
we shall deduce Table 1 (the columns on  $X, \, X^{\mu_p}$)
in an approach different from \S 3.
We shall use the following information:

\smallskip\noindent
(1) $X$  is a del Pezzo surface (Lemma 1.7).
The minimality and rank assumption of the pair imply that
$\mu_p$  acts faithfully on the sublattice 
$M = K_X^\perp$ of $Pic \, X$, and also  $K_X^2 \le 5$,
noting that there are exactly 3 (resp. 6)  $(-1)$-curves on  $X$
when  $K_X^2 = 7$ (resp. $K_X^2 = 6$).

\smallskip\noindent
(2) The lattice $M$ is isomorphic to the root lattice $E_n$,
where $E_5 = D_5, \, E_4 = A_4$.  Here $n = 9-K_X^2 \ge 4$.

\smallskip\noindent
(3) The image $g^*$ of $g$ in $O(E_n)$ belongs to the Weyl group $W(E_n)$.

\smallskip\noindent
(4) All conjugacy classes in Weyl groups are known; see the tables in ATLAS
of finite groups or [Car].

\smallskip\noindent
(5)  $E_n$ embeds naturally into $E_{n+1}$, corresponding to the
natural embeddings of the Dynkin diagrams.  Any ''old'' conjugacy class
in $W(E_n)$ coming from $W(E_{n-1})$, leaves a disjoint union of
$(-1)$-curves invariant and then the pair cannot be minimal
[Man2, Th 6.3].

\smallskip\noindent
(6) Denote by  $C$  the unique, if exists, 
(smooth) irreducible curve in  $X^{\mu_p}$  and write  $C = -m K_X$ 
(Lemma 1.7), where  $2g(C) - 2 = m(m-1)K_X^2$;
we put  $m = 0$  if  $X^{\mu_p}$  is a finite set.
The topological and holomorphic Lefschetz fixed point formulae
in Lemma 1.6 give ($a_3 = 1/4, \, a_4 = 1/2$  when  $p = 5$):  
$$\gather 9 - K_X^2 = (p-1)[3-\sum_{i=1}^{p-1}k_i +m(m-1)K_X^2], \\
1 = \frac{m K_X^2}{12}[(p-2)m + 3] + k_1 \frac{5-p}{12} 
+ k_2 \frac{11-p}{24} + k_3 a_3 + k_4 a_4 + \cdots. \endgather $$

\smallskip\noindent
(7) We have  $(Pic \, X)^{\mu_p} = \bold Z K_X$ (Lemma 1.7).

\par
Now we are in business.

\medskip
Step 1. We know that only $p = 2,3,5,7$ can divide $\#W(E_n)$.

\medskip
Step 2. $p = 7$ can divide only $\#W(E_7), \, \#W(E_8)$. The conjugacy
class  of  $g$  in  $W(E_8)$  is coming from the subgroup  $W(E_7)$
since there is only one each for $n = 7$ and $n = 8$.
So  $n = 7$  and  $K_X^2 = 2$  by the minimality of the pair.
The number of unordered sets of 7 disjoint $(-1)$-curves (an Aronhold set)  
on a degree 2 del Pezzo surface  $X$  is equal to
$\#Sp(6,F_2)/7!$ (see for example [DO, page 167]). Since
the number $36\times 8$ is congruent to 1 mod 7, there is
a $g$-invariant Aronhold set, a contradiction to the minimality
of the pair.

\medskip
Step 3. Assume $p = 5$. It is a new conjugacy class for $n = 4$ $(W(E_4)
= S_5$) and for $n = 8$ (for $n\le 7$ there is only one conjugacy class,
so it is always old).  If  $n = 4$, then  $X$  is a del Pezzo surface
of degree 5.  Hence  $(X, \mu_5)$  fits the last row of Table 1
(Lemma 2.13). 

\medskip
Step 4. Assume $p = 5$ and $n = 8$.  Then  $K_X^2 = 1$.
The formulae in (6) above imply
$(m; k_1, \dots, k_4) = (1;0,0,0,1)$.
So  $(X, \mu_5)$  fits the rows 5 and 6 of Table 1.

\medskip
Step 5. Assume $p = 3$. There is a new conjugacy class of order 3 for
every  $n = 3, 6, 8$. 
If $n = 6$, there is only one new conjugacy class of
order 3; it is  $c_{11}$  in [Man2, Table 1, p. 176];
in Carter's classification it corresponds to the graph
$3A_2\subset E_6$. Its trace on $E_6$ is equal to $-3$
(see also Lemma 1.6).  Now the formulae in (6) imply
$(m; k_1, k_2) = (1; 0, 0), (0; 0, 0)$.
The second case says that  $X^{\mu_3} = \emptyset$, which
is impossible (Lemma 1.2); the first case is the row 3 of Table 1.

\par
Finally, consider the case $n = 8$. There are 2 new
conjugacy classes of order 3.  In Carter's notation they are  $4A_2$  and
$E_8(a_8)$ (they should generate the same group).
As above the formulae in (6) imply  $(m; k_1, k_2) = (2; 1, 0)$.
This is the row 4 of Table 1.

\par \vskip 1.5pc \noindent
Igor V. Dolgachev
\smallskip\noindent
Department of Mathematics
\smallskip \noindent
University of Michigan
\smallskip\noindent
Ann Arbor, MI 48109, USA
\smallskip\noindent
e-mail: idolga$\@$math.lsa.umich.edu

\par \vskip 1.5pc \noindent
{\bf References}

\par \vskip 0.5pc \noindent
[Alex] V. A. Alexeev, Theorems about good divisors on log Fano varieties 
(case of index $r>n-2$), Algebraic geometry (Chicago, IL, 1989),
Lecture Notes in Math. {\bf 1479} (1991), 1--9, Springer, Berlin.

\par \noindent
[Am] F. Ambro, Ladders on Fano varieties, alg-geom/{\bf 9710005}.

\par \noindent
[Art] M. Artin, On isolated rational singularities of surfaces, 
Amer. J. Math. {\bf 88} (1966), 129--136. 

\par \noindent
[ASI, III] M. F. Atiyah and I. M. Singer, The index of elliptic operators, I; III,
Ann. of Math. {\bf 87} (1968), 484--530; 546--604.

\par \noindent
[ASII] M. F. Atiyah and G. B. Segal, The index of elliptic operators, II, 
Ann. of Math. {\bf 87} (1968), 531--545.

\par \noindent
[BPV] W. Barth, C. Peters and A. van de Ven,
Compact complex surfaces, Ergebnisse der Math. und ihrer
Grenzgebiete (3) [Results in Math. and Related Areas (3)], {\bf 4} (1984),
Springer-Verlag, Berlin-New York.

\par \noindent
[BB] L. Bayle, and A. Beauville, Birational involutions of  
$\bold P^2$, math.AG/{\bf 9907028}.

\par \noindent
[Bom] E. Bombieri, Canonical models of surfaces of general type,
Inst. Hautes Etudes Sci. Publ. Math. {\bf 42} (1973), 171--219.

\par \noindent
[BBD] D. Bindschadler, L. Brenton, D. Drucker, 
Rational mappings of del Pezzo surfaces and singular
compactifications of two-dimensional affine varieties, Tohoku Math. J. 
{\bf 36} (1984), 591--609.

\par \noindent
[Bri] E. Brieskorn, Rationale Singularitaten komplexer Flachen, Invent. Math. 
{\bf 4} (1967/1968), 336--358.

\par \noindent
[Car] R.W.Carter, Conjugacy classes in the Weyl group, 
Compos.Math. {\bf 25} (1972), 1-59. 

\par \noindent
[Dem] M. Demazure, Seminaire sur les Singularites des Surfaces, 
Held at the Centre de Math. de l'Ecole Poly., Palaiseau, 1976--1977. 
Lecture Notes in Math., {\bf 777} (1980), Springer, Berlin.

\par \noindent
[DO] I. Dolgachev and D. Ortland, Point sets in projective spaces and 
theta functions, Asterisque {\bf 165} (1988), 210 pp. 

\par \noindent
[FKL] A. Fujiki, R. Kobayashi and S. Lu, 
On the fundamental group of certain open normal surfaces,
Saitama Math. J. {\bf 11} (1993), 15--20.

\par \noindent
[Fuj] T. Fujita, Classification theories of polarized varieties, 
London Mathematical Society Lecture Note Series, {\bf 155} (1990),
Cambridge University Press, Cambridge.

\par \noindent
[Fur] M. Furushima, Singular del Pezzo surfaces and analytic compactifications 
of $3$-dimensional complex affine space ${\bold C}^3$, Nagoya Math. J. 
{\bf 104} (1986), 1--28. 

\par \noindent
[Giz] M. H. Gizatullin, Rational $G$-surfaces,
Math. USSR Izv. {\bf 16} (1981), 103--134.

\par \noindent
[GZ1, 2] R. V. Gurjar and D. -Q. Zhang, 
$\pi\sb 1$ of smooth points of a log del Pezzo surface is finite. I; II, 
J. Math. Sci. Univ. Tokyo {\bf 1} (1994), 137--180; 
{\bf 2} (1995), 165--196.

\par \noindent
[GZ3] R. V. Gurjar and D. -Q. Zhang, On the fundamental groups 
of some open rational surfaces, Math. Ann. {\bf 306} (1996), 15--30.

\par \noindent
[Hart] R. Hartshorne, Algebraic geometry, Graduate Texts in Mathematics, 
{\bf 52} (1977), Springer-Verlag, New York-Heidelberg.

\par \noindent
[Ho1] T. Hosoh, Automorphism groups of quartic del Pezzo surfaces,
J. of Alg. {\bf 185} (1996), 374--389.

\par \noindent
[Ho2] T. Hosoh, Automorphism groups of cubic surfaces,
J. of Alg. {\bf 192} (1997), 651--677.

\par \noindent
[Isk] V. A. Iskovskih, Minimal models of rational surfaces over 
arbitrary fields, Math. USSR Izv. {\bf 14} (1981), 17--39.

\par \noindent
[Kan] S. Kantor, Theorie der endlichen Gruppen von eindeutigen Transformationen
in der Ebene, Berlin: Mayer $\&$ Muller, 1895.

\par \noindent
[Kaw] Y. Kawamata, The cone of curves of algebraic varieties, Ann. of Math. 
{\bf 119} (1984), 603--633.

\par \noindent
[KO] S. Kobayashi and T. Ochiai, Characterizations of complex projective 
spaces and hyperquadrics, J. Math. Kyoto Univ. {\bf 13} (1973), 31--47. 

\par \noindent
[Ko] J. Kollar, Shafarevich maps and plurigenera of algebraic varieties,
Invent. Math. {\bf 113} (1993), 177--215.

\par \noindent
[Koit] M. Koitabashi, Automorphism groups of generic rational surfaces, 
J. of Alg. {\bf 116} (1988), 130--142.

\par \noindent
[Man 1] Yu. I. Manin, Rational surfaces over perfect fields, II, Math. USSR Sb.
{\bf 1} (1967), 141--168.

\par \noindent
[Man 2] Yu. I. Manin, Cubic forms : Algebra, geometry, arithmetic,
2nd ed, North-Holland Math. Library, {\bf 4} (1986),
North-Holland Publ. Co., Amsterdam-New York.

\par \noindent
[MS] S. Keel and J. McKernan, Rational curves on quasi-projective surfaces, 
Mem. Amer. Math. Soc. {\bf 140} (1999).

\par \noindent
[MP] R. Miranda and U. Persson, On extremal rational elliptic surfaces, 
Math. Z. {\bf 193} (1986), 537--558. 

\par \noindent
[MZ] M. Miyanishi and D. -Q. Zhang, Gorenstein log del Pezzo surfaces 
of rank one, J. of Algebra {\bf 118} (1988), 63--84.

\par \noindent
[Mor] S. Mori, Threefolds whose canonical bundles are not numerically effective,
Ann. of Math. {\bf 116} (1982), 133--176. 

\par \noindent
[Mum] D. Mumford, The topology of normal singularities of an algebraic surface 
and a criterion for simplicity, Inst. Hautes Etudes Sci. Publ. Math. 
{\bf 9} (1961), 5--22. 

\par \noindent
[Nag] M. Nagata, On rational surfaces, I, Mem. Coll. Sci. Univ. Kyoto, 
{\bf 32} (1960), 351-370.

\par \noindent
[Re1] M. Reid, Campedelli versus Godeaux, in : Problems in the theory of 
surfaces and their classification (Cortona, 1988),
309--365, Sympos. Math., {\bf XXXII} (1991), Academic Press, London.

\par \noindent
[Re2] M. Reid, Chapters on algebraic surfaces, in:  Complex algebraic geometry 
(Park City, UT, 1993), 3-159, IAS/Park City Math. Ser. {\bf 3}(1997), 
Amer.Math.Soc., Providence, RI.

\par \noindent
[Seg] B. Segre, The Non-singular cubic surfaces, 
Oxford University Press, Oxford, 1942. 

\par \noindent
[Sh] T. Shioda, On the Mordell-Weil lattices, Comment. Math. Univ. St. Paul. 
{\bf 39} (1990), 211--240.

\par \noindent
[Ura] T. Urabe, On singularities on degenerate del Pezzo surfaces of degree $1,$ $2$,
Singularities, Part 2 (Arcata, Calif., 1981), 
Proc. Sympos. Pure Math., {\bf 40} (1983), 587--591, Amer. Math. Soc., 
Providence, R.I.

\par \noindent
[Zh1] D. -Q. Zhang, On Iitaka surfaces, Osaka J. Math. {\bf 24} (1987), 417--460. 

\par \noindent
[Zh2] D. -Q. Zhang, Logarithmic del Pezzo surfaces with rational double 
and triple singular points. Tohoku Math. J. {\bf 41} (1989), 399--452. 

\par \vskip 1.5pc \noindent
D. -Q. Zhang
\smallskip \noindent
Department of Mathematics
\smallskip \noindent
National University of Singapore
\smallskip\noindent
2 Science Drive 2, Singapore 117543
\smallskip\noindent
Republic of Singapore
\smallskip\noindent
e-mail: matzdq$\@$math.nus.edu.sg

$$\text{\bf Table 1}$$

\par \noindent
No.1.  $p \ge 3, \,\, X = {\bold P}^2, \,\, \mu_p = \langle diag[1,1,\zeta_p] \rangle, \,\,
X^{\mu_p} = \{Z = 0\} \cup \{[0,0,1]\}$,
\newline
$Y = X/\mu_p = \overline{{\bold F}}_p, \,\,$ Sing $Y = \frac{1}{p}(1,1), \,\, r = r(Y) = (p+2)/p, \,\,$
details: Ex 2.1a.

\par \vskip 1pc \noindent
No.2.  $p \ge 3, \,\, X = {\bold P}^2, \,\, \mu_p = \langle diag[1,\zeta_p,\zeta_p^v] \rangle
(2 \le v < p), \,\,$
$X^{\mu_p} =$ 
\newline
$\{[1,0,0], [0,1,0], [0,0,1]\}, \,\,
\pi_1(Y^0) = \mu_p, \,\,$ Sing $Y = \frac{1}{p}(1,v), \frac{1}{p}(1,p+1-v),$ 
\newline
$\frac{1}{p}(1,(v-1)/v), 
\,\, r = 3/p \, (p \ge 3)$ or $r = 1/p \, (p \ge 5), \,\,$
details: Ex 2.1b.

\par \vskip 1pc \noindent
No.3. $p = 3, \,\, X$ is cubic del Pezzo, $X^{\mu_p}$ is smooth in $|-K_X|, \,\, Y = {\bold P}^2, \,\,$
\newline
Sing$Y = \emptyset, \,\, r = 3, \,\,$ details: Ex 2.2.

\par \vskip 1pc \noindent
No.4. $p = 3, \,\, X$ is deg 1 del Pezzo, $|-K_X|$ has 6 cuspidal members, 
\newline
$X^{\mu_p}$ is 
a union of a point and a smooth curve of genus 2, $\,\, Y = \overline{{\bold F}}_3, \,\,$
\newline
Sing$Y = \frac{1}{3}(1,1), \,\, r = 5/3, \,\,$ details: Ex 2.3.

\par \vskip 1pc \noindent
No.5. $p = 5, \,\, X$ is deg 1 del Pezzo, $|-K_X|$ has 1 cuspidal and 10 nodal members, $\,\, X^{\mu_p}$ is 
a union of a point and a smooth member in $|-K_X|, \,\, K_Y^2 = 5, \,\, \pi_1(Y^0) = (1), \,\,$
\newline
Sing$Y = \frac{1}{5}(1,4) , \,\, r = 1 \,\,$ details: Ex 2.4, Lemma 2.12.

\par \vskip 1pc \noindent
No.6. $p = 5, \,\, X$ is deg 1 del Pezzo, $|-K_X|$ has 6 cuspidal members, $\,\, X^{\mu_p}$ is 
a union of a point and a smooth member in $|-K_X|, \,\, K_Y^2 = 5, \,\, \pi_1(Y^0) = (1), \,\,$
\newline
Sing$Y = \frac{1}{5}(1,4) , \,\, r = 1 \,\,$ details: Ex 2.4, Lemma 2.12.

\par \vskip 1pc \noindent
No.7. $p = 5, \,\, X$ is the deg 5 del Pezzo, $X^{\mu_p}$ is 
a 2-point set, $\,\, K_Y^2 = 1, \,\, \pi_1(Y^0) = \mu_5, \,\,$
\newline
Sing$Y = \frac{1}{5}(1,4) , \frac{1}{5}(1,4), \,\, r = 1 \,\,$ details: Ex 2.5, Lemma 2.13.

$$\text{\bf Table 2}$$

\par \noindent
No.1. $G = {\bold Z}/(3), \,\, X = {\bold P}^2, \,\, G = \langle diag[1, \zeta_3, \zeta_3^2] \rangle, \,\,
K_Y^2 = 3, \,\,$ 
\newline
Sing$Y = \frac{1}{3}(1,2), \frac{1}{3}(1,2), \frac{1}{3}(1,2), \,\,$
details: Ex 2.1b.

\par \vskip 1pc \noindent
No.2. $G = {\bold Z}/(4), \,\, X = {\bold P}^1 \times {\bold P}^1, \,\, 
X^G$ is a 4-point set, $\,\,K_Y^2 = 2, \,\,$ 
\newline
Sing$Y = \frac{1}{2}(1,1), \frac{1}{4}(1,3), 
\frac{1}{4}(1,3), \,\,$
details: Ex 2.11, Lemma 2.15.

\par \vskip 1pc \noindent
No.3. $G = {\bold Z}/(5), \,\, X$ is the deg 5 del Pezzo, $\,\, 
X^G$ is a 2-point set, $\,\,K_Y^2 = 1, \,\,$ Sing$Y = \frac{1}{5}(1,4), \frac{1}{5}(1,4), \,\,$
details: Ex 2.5, Lemma 2.13.

\par \vskip 1pc \noindent
No.4. $G = {\bold Z}/(6), \,\, X$ is the deg 6 del Pezzo, $\,\, 
X^G$ is a 6-point set, $\,\,K_Y^2 = 1, \,\,$ 
\newline
Sing$Y = \frac{1}{2}(1,1), \frac{1}{3}(1,2), 
\frac{1}{6}(1,5), \,\,$
details: Ex 2.9, Lemma 2.14.

\par \vskip 1pc \noindent
No.5. $G = {\bold Z}/(3) \oplus {\bold Z}/(3), \,\, X = {\bold P}^2, \,\, 
X^G$ is a 12-point set, $\,\,
G = \langle g_1, g_2 \rangle$,
\newline 
$g_1 = diag[1, \zeta_3, \zeta_3^2], \,\,
g_2 = (a_{ij}) (a_{13}a_{21}a_{32} = 1$, other $a_{ij} = 0), \,\,
K_Y^2 = 1, \,\,$ 
\newline
Sing$Y = \frac{1}{3}(1,2), \frac{1}{3}(1,2), 
\frac{1}{3}(1,2), \frac{1}{3}(1,2), \,\,$
details: Ex 2.10.

\par \vskip 1pc \noindent
No.6. $G = {\bold Z}/(4) \oplus {\bold Z}/(2), \,\, X = {\bold P}^1 \times {\bold P}^1$,
$X^G$ is a 12-point set, $\,\,K_Y^2 = 1, \,\,$ 
\newline
Sing$Y = \frac{1}{2}(1,1), \frac{1}{2}(1,1), 
\frac{1}{4}(1,3), \frac{1}{4}(1,3), \,\,$
details: Ex 2.11, Lemma 2.15.

\enddocument

\documentstyle[11pt]{article}
\textwidth 6in
\textheight 9in
\oddsidemargin 0.3cm
\evensidemargin 0.3cm
\topmargin -1.5cm
\begin{document}

$${\bf Table 1}$$
\begin{tabular}{|l|l|l|l|l|l|r|} \hline
$p$                                          &$X$                                      
&$X^{\mu_p}$; \, $\mu_p = \langle g \rangle$ &$Y=X/\mu_p$; $\pi_1$   &Sing $Y$               
&$r = r(Y)$                                  &details \\ \hline \hline

                                           &
&{\rm line} \, $\{Z=0\}$                   & & & &   \\
$p \geq 3$                                 &${\bf P^2}$
& \& \, {\rm point} $\,$ [0,0,1]              &$\overline{\bf F}_p$     &$\frac{1}{p}(1,1)$
&$(p+2)/p$                                 & Ex2.1a \\
                                           &
&$g = {\rm diag} [1, 1, \zeta_p]$          & & & & \\ \hline

& &$[1,0,0], [0,1,0]$                      &                         &$\frac{1}{p}(1,v)$
&$r=\frac{3}{p}$($p\ge3$)                              & \\
$p \geq 3$                                 &${\bf P^2}$
&[0,0,1]; $2\le v< p$                      &$\pi_1(Y^0) = \mu_p$     &$\frac{1}{p}(1,p+1-v)$  
&or                                        &Ex2.1b \\
                                           &
&$g = {\rm diag} [1, \zeta_p, \zeta_p^v]$  &  
&$\frac{1}{p}(1,(v-1)/v)$                  &$r=\frac{1}{p}$($p\ge5$) & \\ \hline

$3$                                        &cubic delPezzo                                    
&smooth $\in |-K_X|$                       &${\bf P^2}$              &$\emptyset$               
&$3$                                       &Ex 2.2    \\ \hline

                                           &deg 1 delPezzo
&smooth                                    & & & &    \\
$3$                                        &$|-K_X|$ has             &of genus 2
&$\overline{\bf F}_3$ &$\frac{1}{3}(1,1)$  &$5/3$                    & Ex 2.3 \\ 
                                           &6 cuspidals              &$\& \,$ {\rm a point} 
& & & & \\ \hline

                                           &deg 1 delPezzo           & &$K_Y^2=5$ & & & \\
$5$                                        &$|-K_X|$ has             &smooth $\in |-K_X|$     
&                                          &$\frac{1}{5}(1,4)$ 
&$1$                                       &Ex 2.4 \\
                                           &10 nodals \&             & \& a point
&$\pi_1(Y^0)=(1)$ & & &2.12 \\
                                           &1 cuspidal               & & & & & \\ \hline

                                           &deg 1 delPezzo           & &$K_Y^2=5$ & & &Ex 2.4 \\
$5$                                        &$|-K_X|$ has             &smooth $\in |-K_X|$     
&$\pi_1(Y^0)=(1)$                          &$\frac{1}{5}(1,4)$ 
&$1$                                       &2.12 \\
                                           &6 cuspidals              &$\&$ a point & & & & \\ \hline

& & &$K_Y^2=1$                             &$\frac{1}{5}(1,4)$       & &Ex 2.5 \\
$5$                                        &deg 5 delPezzo
&2 points                                  &$\pi_1(Y^0)=\mu_5$       &$\frac{1}{5}(1,4)$ 
&$1$                                       &2.13 \\ \hline \hline
\end{tabular}
\end{document}